\documentclass[11pt]{amsart}
\usepackage{mathrsfs}
\usepackage{amssymb}

\usepackage{titletoc}
\usepackage{amscd}
\usepackage{amsmath, amssymb}
\usepackage{amsfonts}
\usepackage[colorlinks,linkcolor=blue, pdfstartview=FitH]{hyperref}

\numberwithin{equation}{section}

\theoremstyle{plain}
\newtheorem{thm}{Theorem}[section]
\newtheorem{theorem}[thm]{Theorem}
\newtheorem{lemma}[thm]{Lemma}
\newtheorem{corollary}[thm]{Corollary}
\newtheorem{proposition}[thm]{Proposition}

\theoremstyle{definition}

\newtheorem{problem}[thm]{Problem}
\newtheorem{remark}[thm]{Remark}

\newtheorem{definition}[thm]{Definition}

\newtheorem{example}[thm]{Example}

\newtheorem{defn-thm}[thm]{Definition-Theorem}

\newcommand{\sO}{{\mathcal O}}

\newcommand{\C}{{\mathbb C}}

\renewcommand{\P}{{\mathbb P}}

\newcommand{\R}{{\mathbb R}}

\newcommand{\T}{{\mathbb T}}

\renewcommand{\S}{{\mathbb S}}

\newcommand{\Z}{{\mathbb Z}}

\newcommand{\tr}{{ tr}}

\newcommand{\qtq}[1]{\quad\mbox{#1}\quad}
\newcommand{\bp}{\bar{\partial}}
\newcommand{\Om}{\Omega}

\newcommand{\ds}{\oplus}

\newcommand{\ts}{\otimes}

\newcommand{\btheorem}{\begin{theorem}}
\newcommand{\etheorem}{\end{theorem}}
\newcommand{\bproposition}{\begin{proposition}}
\newcommand{\eproposition}{\end{proposition}}
\newcommand{\bdefinition}{\begin{definition}}
\newcommand{\edefinition}{\end{definition}}
\newcommand{\bcorollary}{\begin{corollary}}
\newcommand{\ecorollary}{\end{corollary}}
\newcommand{\bproof}{\begin{proof}}
\newcommand{\eproof}{\end{proof}}
\newcommand{\bremark}{\begin{remark}}
\newcommand{\eremark}{\end{remark}}
\newcommand{\eexample}{\end{example}}
\newcommand{\bexample}{\begin{example}}
\newcommand{\la}{\langle}
\newcommand{\elemma}{\end{lemma}}
\newcommand{\blemma}{\begin{lemma}}
\newcommand{\ra}{\rangle}
\newcommand{\sq}{\sqrt{-1}}

\newcommand{\p}{\partial}

\renewcommand{\bar}{\overline}

\renewcommand{\phi}{\varphi}

\newcommand{\ee}{\end{eqnarray*}}
\newcommand{\be}{\begin{eqnarray*}}

\newcommand{\beq}{\begin{equation}}
\newcommand{\eeq}{\end{equation}}

\newcommand{\bd}{\begin{enumerate}}
\newcommand{\ed}{\end{enumerate}}

\renewcommand{\hat}{\widehat}
\renewcommand{\tilde}{\widetilde}

\renewcommand{\>}{\rightarrow}

\newcommand{\fR}{\mathfrak R}

\begin{document}
\title{Ricci  curvatures on Hermitian manifolds}

\makeatletter
\let\uppercasenonmath\@gobble
\let\MakeUppercase\relax
\let\scshape\relax
\makeatother
\author{Kefeng Liu, Xiaokui Yang}
\date{}

\address{Kefeng Liu, Department of Mathematics, UCLA, Los Angeles, CA}
\email{liu@math.ucla.edu}

\address{{ Xiaokui Yang, Morningside Center of Mathematics, Institute of
        Mathematics, Hua Loo-Keng Key Laboratory of Mathematics,
        Academy of Mathematics and Systems Science,
        Chinese Academy of Sciences, Beijing, 100190, China.}}
\email{\href{mailto:xkyang@amss.ac.cn}{\texttt{xkyang@amss.ac.cn}}}

\maketitle

\begin{abstract}
In this paper, we introduce the first Aeppli-Chern class for complex
manifolds and show that the $(1,1)$- component of the curvature
$2$-form of the Levi-Civita connection on the anti-canonical line
bundle represents this class. We systematically investigate the
relationship between a variety of Ricci curvatures on Hermitian
manifolds and the background Riemannian manifolds. Moreover, we
study non-K\"ahler Calabi-Yau manifolds by using the first
Aeppli-Chern class and the Levi-Civita Ricci-flat metrics. In
particular, we construct explicit Levi-Civita Ricci-flat metrics on
Hopf manifolds $\S^{2n-1}\times \S^1$. We also construct a smooth
family of Gauduchon metrics on a compact Hermitian manifolds such
that the metrics are in the same first Aeppli-Chern class, and their
first Chern-Ricci curvatures are the same and nonnegative, but their
Riemannian scalar curvatures are constant and vary smoothly between
 negative infinity and a positive number. In particular, it shows
 that
 Hermitian manifolds with nonnegative first Chern class can admit  Hermitian metrics with strictly negative Riemannian scalar
curvature.

\end{abstract}

{{\scriptsize{\setcounter{tocdepth}{2} \tableofcontents
\dottedcontents{section}[0.8cm]{}{1.8em}{5pt}
\dottedcontents{subsection}[1.70cm]{}{2.7em}{5pt}
}}



\section{Introduction}

%
%
%

 Let $(M,h)$ be a Hermitian manifold and $g$ the background Riemannian metric. It is well-known that, when $(M,h)$ is not K\"ahler, the complexification of the real curvature tensor $R$
 is extremely complicated. Moreover,
 on the Hermitian holomorphic vector bundle $(T^{1,0}M,h)$, there are two
 typical connections: the (induced) Levi-Civita connection and the Chern
 connection.
 The curvature tensors of them are denoted by $\mathfrak R$ and
 $\Theta$ respectively. It is known that the complexified
 Riemannian curvature $R$, the Hermitian Levi-Civita curvature $\mathfrak R$ and the Chern
 curvature $\Theta$ are mutually different. It is also obvious
 that $R$ is closely related to the Riemannian geometry of $M$, and
 $\Theta$ can characterize many complex geometric properties of $M$ whereas  $\mathfrak
 R$ can be viewed as a bridge between $R$ and $\Theta$, i.e. a bridge between Riemannian geometry and Hermitian geometry.

   Let
 $\{z^i\}_{i=1}^n$ be the local holomorphic coordinates centered at
 a point $p\in M$.  We can compare the curvature tensors of $\mathfrak R$, $\Theta$ when restricted on the space $\Gamma(M, \Lambda^{1,1}T^*M\ts
 End(T^{1,0}M))$. That is, we can find relations between  $\mathfrak R_{i\bar j k\bar \ell}$ and $\Theta_{i\bar j k\bar
 \ell}$. We denoted
 by $$\mathfrak R^{(1)}=\sq \mathfrak R^{(1)}_{i\bar j} dz^i\wedge d\bar z^j\qtq{with} \mathfrak R^{(1)}_{i\bar j}=h^{k\bar\ell} \mathfrak R_{i\bar j k\bar \ell},$$
$$\mathfrak R^{(2)}=\sq \mathfrak R^{(2)}_{i\bar j} dz^i\wedge d\bar z^j\qtq{with} \mathfrak R^{(2)}_{i\bar j}=h^{k\bar\ell} \mathfrak R_{k\bar \ell i\bar j },$$
$$\mathfrak R^{(3)}=\sq \mathfrak R^{(3)}_{i\bar j} dz^i\wedge d\bar z^j\qtq{with} \mathfrak R^{(3)}_{i\bar j}=h^{k\bar\ell} \mathfrak R_{i\bar \ell k \bar j},$$
and
$$\mathfrak R^{(4)}=\sq \mathfrak R^{(4)}_{i\bar j} dz^i\wedge d\bar z^j\qtq{with} \mathfrak R^{(4)}_{i\bar j}=h^{k\bar\ell} \mathfrak R_{k\bar j i\bar \ell }.$$
$\mathfrak R^{(1)}$ and  $\mathfrak R^{(2)}$ are called the
\emph{first Levi-Civita Ricci curvature} and the \emph{second
Levi-Civita Ricci curvature} of $(T^{1,0}M, h)$ respectively;
$\mathfrak R^{(3)}$and $\mathfrak R^{(4)}$ are the corresponding
third and fourth Levi-Civita Ricci curvatures. Similarly, we can
define the first Chern-Ricci curvature $\Theta^{(1)}$, the second
Chern-Ricci curvature $\Theta^{(2)}$, the third and fourth
Chern-Ricci curvatures $\Theta^{(3)}$ and $\Theta^{(4)}$
respectively. As shown in \cite{LY12}, $\mathfrak R^{(2)}$ and
$\Theta^{(2)}$ are closely related to the geometry of $M$, for
example, we can use them to study the cohomology groups and
plurigenera of compact Hermitian manifolds. On the other hand, it is
well-known that $\Theta^{(1)}$ represents the first Chern class $
c_1(M)\in H^{1,1}_{\bp}(M)$. However, in general, the first
Levi-Civita Ricci form $\mathfrak R^{(1)}$ is not $d$-closed, and so
it can not represent a class in $H^{1,1}_{\bp}(M)$.\\

  We introduce two cohomology groups to study the geometry  of
  compact complex (especially, non-K\"ahler) manifolds, the Bott-Chern
  cohomology and the Aeppli cohomology:
  $$ H^{p,q}_{BC}(M):= \frac{\text{Ker} d \cap \Om^{p,q}(M)}{\text{Im} \p\bp \cap \Om^{p,q}(M)}\qtq{and} H^{p,q}_{A}(M):=\frac{\text{Ker} \p\bp \cap \Om^{p,q}(M)}{\text{Im} \p \cap \Om^{p,q}(M)+ \text{Im}\bp \cap
  \Om^{p,q}(M)}.$$
Suppose $\alpha$ is a $d$-closed $(p,q)$-form. We denote by
$[\alpha]_{BC}$ and $[\alpha]_A$, the corresponding classes in
$H^{p,q}_{BC}(M)$ and $H^{p,q}_A(M)$ respectively.  Let
$\mathrm{Pic}(M)$ be the set of holomorphic line bundles over $M$.
As similar as the first Chern class map
$c_1:\mathrm{Pic}(M)\>H^{1,1}_{\bp}(M)$, there is a \textbf{first
Aeppli-Chern class} map \beq c_1^{{AC}}:\mathrm{Pic}(M)\>
H^{1,1}_{{A}}(M),\eeq which can be described as follows.

\bdefinition Let $L\to M$ be a holomorphic line bundle over $M$. The
first Aeppli-Chern class is defined as \beq
c^{AC}_1(L)=\left[-\sq\p\bp\log h\right]_A \in H^{1,1}_{A}(M)\eeq
where $h$ is an arbitrary smooth Hermitian metric on $L$.  Note that
$-\sq\p\bp\log h$ is the (local) curvature form $\Theta_h$ of the
Hermitian line bundle $(L,h)$.  If we choose a different metric
$h'$, then $\Theta_{h'}-\Theta_h=\sq\p\bp\log
\left(\frac{h}{h'}\right)$ is globally $\p\bp$-exact.  Hence
$c_1^{{AC}}(L)$  is well-defined  in $H^{1,1}_{{A}}(M)$ and it is
independent of the metric $h$.  For a complex manifold $M$,
$c_1^{{AC}}(M)$ is defined to be $c_1^{{AC}}(K^{-1}_M)$ where
$K_M^{-1}$ is the anti-canonical line bundle $\wedge^n T^{1,0}M$.
\edefinition

\noindent Note that, for a Hermitian line bundle $(L,h)$, the
classes $c_1(L)$, $c^{BC}_1(L)$ and $c_1^{AC}(L)$ have the same
$(1,1)$-form
representative $\Theta^h=-\sq\p\bp\log h$ (in different classes).\\

 On K\"ahler manifolds or Hermitian manifolds the first Chern classes and first Bott-Chern classes are
 well-studied in the literatures by using the Chern connection. In particular, the related
 Monge-Amp\`ere type equations are extensively investigated since the celebrated
 work of Yau.  However, the geometry
 of the Levi-Civita connection is not well understood in complex
 geometry although it has rich real Riemannian geometry
 structures.

\subsection{The complex geometry of the Levi-Civita connection}

In the first part of this paper, we study (non-K\"ahler) Hermitian
manifolds by using
  Levi-Civita Ricci
forms  and Aeppli-Chern class $c_1^{AC}(M)$. At first, we establish
the following  result which is analogous to the classical result
that  the (first) Chern-Ricci curvature $\Theta^{(1)}$ represents
the first Chern class $c_1(M)$.
 \btheorem\label{main-1} On a compact Hermitian manifold $(M,\omega)$, the first Levi-Civita Ricci
form $\fR^{(1)}$ represents the first Aeppli-Chern class $
c_1^{AC}(M)$ in $H^{1,1}_{A}(M)$. More precisely, \beq
\mathfrak{R}^{(1)}=\Theta^{(1)}-\frac{1}{2}(\p\p^*\omega+\bp\bp^*\omega)\label{00}.\eeq
In particular, we obtain
\begin{enumerate}
\item $\mathfrak R^{(1)}=\Theta^{(1)}$ if and only if $d^*\omega=0$,
i.e. $(M,\omega)$ is a balanced manifold;

\item if
$\bp\p^*\omega=0$, then $\mathfrak{R}^{(1)}$ represents the real
first Chern class $c_1(M)\in H^{2}_{dR}(M)$, i.e.
$c_1(M)=c_1^{AC}(M)$ in $H^{2}_{dR}(M)$.
\end{enumerate} \etheorem

\noindent We also show that, on complex manifolds supporting
$\p\bp$-lemma (e.g. manifolds in Fujiki class $\mathscr C$ and in
particular, Moishezon manifolds), the converse statement of $(2)$
holds.\vspace{0.2cm}

 There is  an important class of manifolds, so called
Calabi-Yau manifolds which are  extensively studied by
mathematicians and also physicians. In this paper, a Calabi-Yau
manifold is a complex manifold with $c_1(M)=0$ and we will focus on
non-K\"ahler Calabi-Yau manifolds.
 There are many fundamental results on non-K\"ahler
Hermitian manifolds  with vanishing first Bott-Chern classes. They
are always characterized by using the first Chern-Ricci curvature
$\Theta^{(1)}$ and  the related Monge-Amp\`ere type equations
(e.g.\cite{FWW10, Gi,TW1,TW3,TW4,TW5,TW6}).
  For more details on this subject, we refer the reader to the nice survey paper
  \cite{T}.\vspace{0.3cm}

 Next, we make the following
 observation:
 \bcorollary\label{c1}  Let $M$ be a complex
manifold. Then
$$c_1^{BC}(M)=0 \Longrightarrow c_1(M)=0 \Longrightarrow c_1^{AC}(M)=0.$$ Moreover, on a complex manifold
satisfying the $\p\bp$-lemma,
$$c_1^{BC}(M)=0\Longleftrightarrow c_1(M)=0\Longleftrightarrow
c_1^{AC}(M)=0.$$
 \ecorollary
\noindent That means, it is very natural to study non-K\"ahler
Calabi-Yau manifolds by using the first Aeppli-Chern class
$c_1^{AC}$ and the first Levi-Civita Ricci curvature $\fR^{(1)}$.

\bdefinition A Hermitian metric $\omega$ on $M$ is called
\emph{Levi-Civita Ricci-flat} if $$\fR^{(1)}(\omega)=0.$$
\edefinition

\noindent It is known that  Hopf manifolds $M=\S^{2n-1}\times \S^1$
are all non-K\"ahler Calabi-Yau manifolds $(n\geq 2)$, i.e.
$c_1(M)=0$.  However,  there is no Chern Ricci-flat Hermitian
metrics on $M$, i.e. there does \textbf{not} exist a Hermitian
metric $\omega$ such that $\Theta^{(1)}(\omega)=0$ since
$c_1^{BC}(M)\neq 0$ (see Remark \ref{remark2} for more details).
 On the contrary, we can construct explicit Levi-Civita Ricci-flat metrics on
 them.

\btheorem\label{example1} Let $M=\S^{2n-1}\times \S^1$ with $n\geq
2$ and $\omega_0$ the canonical metric on $M$.  The perturbed
Hermitian metric $$ \omega=\omega_0-\frac{4}{n}\fR^{(1)}(\omega_0)$$
is Levi-Civita Ricci-flat, i.e. $\fR^{(1)}(\omega)=0$. \etheorem

\noindent
 Theorem \ref{example1} is an explicit example demonstrating that the
 zero first Aeppli-Chern class can imply the existence of  Levi-Civita Ricci-flat Hermitian
 metric. This result is also anologous to the following two
 classical and general results. One is  Yau's celebrated solution to Calabi's
 conjecture
\btheorem[{\cite{Yau}}]\label{yau} Let $(M,\omega)$ be a compact
K\"ahler manifold. If the real $(1,1)$ form $\eta$ represents the
first Chern class $c_1(M)$, then there exists a smooth function
$\phi\in C^\infty(M)$ such that the K\"ahler metric $\tilde
\omega=\omega+\sq\p\bp\phi$  has Ricci curvature $\eta$, i.e. $$
Ric(\tilde\omega)=\eta.$$ \etheorem

\noindent The other one is  Tosatti and Weinkove's Hermitian
analogue of Yau's fundamental result:

\btheorem[{\cite{TW1}}]\label{TW} Let $(M,\omega)$ be a compact
Hermitian manifold. If the real $(1,1)$ form $\eta$ represents the
first Bott-Chern class $c^{BC}_1(M)$, then there exists a smooth
function $\phi\in C^\infty(M)$ such that the Hermitian metric
$\tilde \omega=\omega+\sq\p\bp\phi$  has Chern-Ricci curvature
$\eta$, i.e. $$ \Theta^{(1)}(\tilde\omega)=\eta.$$ \etheorem

As inspired by Theorem \ref{example1} and the fundamental Theorem
\ref{yau} and Theorem \ref{TW}, we propose the following problem.
\begin{problem}\label{question1}  Let $M$ be a compact complex manifold. For a fixed Hermitian metric $\omega_0$ on $M$,
 and a real $(1,1)$-form $\eta$ representing $c^{AC}_1(M)$, does there exist
a $(0,1)$-form $\gamma$ such that the Hermitian metric $
\omega=\omega_0+\p\gamma+\bp\bar\gamma$ satisfies
$$\fR^{(1)}(\omega)=\eta?$$ In particular, if
$c_1^{AC}(M)=0$ (or $c_1(M)=0$), does there exist a $(0,1)$-form
$\gamma$ such that the Hermitian metric $
\omega=\omega_0+\p\gamma+\bp\bar\gamma$ is Levi-Civita Ricci-flat,
i.e. $\fR^{(1)}(\omega)=0$?
\end{problem}
\noindent Although the background PDE is not exactly the same as the
standard Monge-Amp\`ere equations, we hope similar methods could
work.

The first Levi-Civita Ricci curvature $\fR^{(1)}$ is  closely
related to possible symplectic structures on Hermitian manifolds. It
is easy to see that, on a Hermitian manifold $(M,\omega)$,
$\fR^{(1)}$ is the $(1,1)$-component of the curvature $2$-form of
the Levi-Civita connection on $K_M^{-1}$. Hence, if $\fR^{(1)}$ is
 strictly positive, it can induce a
symplectic structure on $M$ (see Theorem \ref{main111} and also
\cite{LZ09}).
 Moreover, the symplectic structures  thus obtained are not
necessarily K\"ahler. \vspace{0.2cm}

As applications of Theorem \ref{main-1}, we can characterize
Hermitian manifolds by using the
 Levi-Civita Ricci curvature.

\btheorem\label{main-100} Let $(M,h)$ be a compact Hermitian
manifold. If the first Levi-Civita Ricci curvature $\fR^{(1)}$ is
quasi-positive, then the top intersection number $c^n_1(M)> 0$.  In
particular, $H^2_{dR}(M)$, $H^{1,1}_{\bp}(M), H^{1,1}_{BC}(M)$ and
$H^{1,1}_A(M)$ are all non-zero. \etheorem

The Levi-Civita Ricci curvature and the Aeppli-Chern class are also
closely related to the algebraic aspects of the anti-canonical line
bundle. For example, on a projective manifold $M$, $c_1(M)$ and
$c^{AC}_1(M)$ are ``numerically" equivalent, i.e. for any
irreducible curve $\gamma$ in $M$,
$$c^{AC}_1(M)\cdot \gamma=c_1(M)\cdot \gamma.$$
 In particular, \bcorollary  \bd\item if  $\fR^{(1)}$ is semi-positive,   the anti-canonical line bundle $K^{-1}_{M}$ is
\emph{nef};\item if  $\fR^{(1)}$ is quasi-positive, then $K_M^{-1}$
is a \emph{big} line bundle.\ed \ecorollary

\noindent It is not hard to see that if the Hermitian manifold
$(M,g)$ has positive constant Riemannian sectional curvature, then
$\fR^{(1)}$ is positive.  On the other hand, since the positivity
condition is an open condition, $\fR^{(1)}$ is still positive in a
small neighborhood of a positive constant sectional curvature
metric. As an application of this  observation,  one can see the
following result of Lebrun which is also observed in \cite{BHL99}
and \cite{Ta06}.

\bcorollary\label{c2} On $\S^6$, there is no orthogonal complex
structure compatible with metrics in  some small neighborhood of the
round metric. \ecorollary

\subsection{Curvature relations on Hermitian manifolds} In the second part of this paper, we investigate the relations between various Ricci curvatures on
Hermitian manifolds. As introduced above, on a  Hermitian manifold
$(M,\omega)$ there are several different types of Ricci curvatures:
\bd
\item the  Levi-Civita Ricci curvatures $\mathfrak{R}^{(1)}$, $\mathfrak{R}^{(2)}$, $\mathfrak{R}^{(3)}$, $\mathfrak{R}^{(4)}$;
\item the Chern Ricci curvatures $\Theta^{(1)}$, $\Theta^{(2)}$, $\Theta^{(3)}$,
$\Theta^{(4)}$;
\item the Hermitian-Ricci curvature $Ric_H=\sq R_{i\bar j} dz^i\wedge d\bar z^j$ where $R_{i\bar j}=h^{k\bar \ell} R_{i\bar j k\bar
\ell}$ (which is equal to $R^{(1)}$ and $R^{(2)}$), the third and
fourth Hermitian-Ricci curvatures $R^{(3)}$ and $R^{(4)}$;
\item the $(1,1)$-component of the complexified
Riemannian Ricci curvature,  $\mathscr Ric$.
 \ed

 \noindent If  $(M,\omega)$ is K\"ahler, all Ricci
curvatures are the same, but it is not true on general Hermitian
manifolds. We shall explore explicit relations between them by using
the Hermitian metric $\omega$ and its torsion $T$. We write them
down with a reference curvature, e.g. $\Theta^{(1)}$, to the
reader's convenience.

\btheorem\label{riccire} Let $(M,\omega)$ be a compact Hermitian
manifold. \bd
\item The  Levi-Civita Ricci curvatures are  $$
\mathfrak{R}^{(1)}=\Theta^{(1)}-\frac{1}{2}\left(\p\p^*\omega+\bp\bp^*\omega\right);$$$$
\mathfrak{R}^{(2)}=\Theta^{(1)}-\frac{1}{2}\left(
\p\p^*\omega+\bp\bp^*\omega\right)-\frac{\sq}{4}T\circ\bar
T+\frac{\sq}{4} T\boxdot \bar T;$$
$$
\mathfrak{R}^{(3)}=\Theta^{(1)}-\frac{1}{2}\left(\sq\Lambda\left(\p\bp\omega\right)+(\p\p^*\omega+\bp\bp^*\omega)\right)+\frac{\sq}{4}
T\boxdot \bar T+\frac{T([\p^*\omega]^{\#})}{4};$$$$
\mathfrak{R}^{(4)}=\Theta^{(1)}-\frac{1}{2}\left(\sq\Lambda\left(\p\bp\omega\right)+(\p\p^*\omega+\bp\bp^*\omega)\right)+\frac{\sq}{4}
T\boxdot \bar T+\frac{\bar{T([\p^*\omega]^{\#})}}{4},$$ \noindent
where $(\p^*\omega)^{\#}$ is the dual vector of  the $(0,1)$-form
$\p^*\omega$.

\item  The  Chern-Ricci curvatures are $$
\Theta^{(2)}=\Theta^{(1)}-\sq\Lambda\left(\p\bp\omega\right)-(\p\p^*\omega+\bp\bp^*\omega)+{\sq}T\boxdot\bar
T;$$$$
\Theta^{(3)}=\Theta^{(1)}-\p\p^*\omega;$$$$\Theta^{(4)}=\Theta^{(1)}-\bp\bp^*\omega.$$

\item The Hermitian-Ricci curvatures are $$  Ric_H=R^{(1)}=R^{(2)}=\Theta^{(1)}-\frac{1}{2}\left(
\p\p^*\omega+\bp\bp^*\omega\right)-\frac{\sq}{4}T\circ\bar T;$$
\begin{eqnarray}
R^{(3)}=R^{(4)}\nonumber&=&\Theta^{(1)}-\frac{1}{2}\left(\sq\Lambda\left(\p\bp\omega\right)+(\p\p^*\omega+\bp\bp^*\omega)\right)\\
\nonumber&&+\frac{\sq}{4} T\boxdot \bar
T+\frac{\bar{T([\p^*\omega]^{\#})}+T([\p^*\omega]^{\#})}{4}.\end{eqnarray}

\item The $(1,1)$-component of the complexified Riemannian Ricci
curvature is \begin{eqnarray} \mathscr
{R}ic\nonumber&=&\Theta^{(1)}-\sq(\Lambda
\p\bp\omega)-\frac{1}{2}(\p\p^*\omega+\bp\bp^*\omega)+\frac{\sq}{4}\left(2T\boxdot
\bar T+T\circ \bar T\right)\\ \nonumber&&+\frac{1}{2}\left(
T([\p^*\omega]^{\#})+\bar{
T([\p^*\omega]^{\#})}\right).\end{eqnarray}\ed

 \etheorem

\noindent From these curvature relations, one can see clearly the
geometry of many Hermitian manifolds with special metrics (e.g.
$d^*\omega=0$, $\p\bp\omega=0$). In particular, these curvature
relations may
 enlighten the study of various Hermitian Ricci flows (e.g.
\cite{ST1,ST2,ST3}, \cite{LY12}, \cite{Gi,TW3,TW4, TWY, Gill1,GS})
by using the well-studied Hamilton's Ricci flow. Of course, it is
also natural to define new Hermitian Ricci flows by certain Ricci
curvatures with significant geometric meanings, for example, the
first Levi-Civita Ricci curvature $\fR^{(1)}$.\vspace{0.3cm}

  As straightforward consequences, from the Ricci curvature relations, we can also obtain relations
between the corresponding scalar curvatures. It is known that the
 positive scalar curvature can characterize the geometry of the
 manifolds. In \cite{Yau2}, Yau proved that, on a compact K\"ahler manifold $(M,\omega)$, if
the total scalar curvature is positive, then all plurigenera
$p_m(M)$ vanish, and so the Kodaira dimension of $M$ is $-\infty$.
Based on Yau's result, Heier-Wong (\cite{HW12}) observed that on a
projective manifold, if the total scalar curvature is positive, the
manifold is uniruled, i.e. it is covered by rational curves. On a
compact Hermitian manifold $M$, Gauduchon showed (\cite{Ga2}) that
if the total Chern scalar curvature of a \emph{Gauduchon metric} is
positive, then $p_m(M)=0$ and $\kappa(M)=-\infty$. On the other
hand, the Riemannian scalar curvature on Riemannian manifolds is
extensively studied. In particular, by Trudinger, Aubin and Schoen's
solution to the Yamabe problem, it is well-known that every
Riemannian metric is conformal to a metric with constant scalar
curvature. To understand the relations between Riemannian geometry
and Hermitian geometry, the following relation
 is
of particular interest.

\bcorollary\label{sphere} On a compact Hermitian manifold
$(M,\omega)$, the Riemannian scalar curvature $s$ and the Chern
scalar curvature $s_C$ are related by \beq
s=2s_C+\left(\la\p\p^*\omega+\bp\bp^*\omega,\omega\ra-2|\p^*\omega|^2\right)-\frac{1}{2}|T|^2.\label{sc}\eeq
Moreover,  according to the different types of Ricci curvatures,
there are  different scalar curvatures and the following statements
are equivalent: \bd

\item $(M,\omega)$ is K\"ahler;

\item $\displaystyle \int s\cdot \omega^n=\int 2s_C\cdot \omega^n$;

\item $\displaystyle \int s_C\cdot \omega^n=\int s_R\cdot \omega^n$;

\item $\displaystyle \int s_C\cdot \omega^n=\int s_{H}\cdot \omega^n$;

\item $\displaystyle \int s_H\cdot \omega^n=\int s_{LC}\cdot
\omega^n$.\ed

 \ecorollary

\noindent A similar formulation as (\ref{sc}) by using ``Lee forms"
is also observed by Gauduchon (\cite{Ga3}). See also \cite{TV81,
Va83, AD99, FG, FPS, AI01, IP01, GI01, AD03,  MS06, IP13} for some
curvature relations on complex surfaces.  For more scalar curvature
relations, see Corollary \ref{A}, Remark \ref{remark3}, Corollary
\ref{B} and Corollary \ref{C}.\vspace{0.2cm}

\subsection{Special metrics on Hermitian manifolds} Finally, we study  special metrics on Hermitian manifolds.
In the following, we give precise examples of Hermitian manifolds on
the relations between Ricci curvatures, Chern scalar curvatures and
Riemannian scalar curvatures.

\btheorem\label{example2} For any $n\geq 2$, there exists an
$n$-dimensional compact complex manifold $X$ with $c_1(X)\geq 0$,
such that $X$ admits three different \textbf{Gauduchon metrics}
$\omega_1,\omega_2$ and $\omega_3$ with the following properties.
\bd
\item $[\omega_1]=[\omega_2]=[\omega_3]\in H^{1,1}_{A}(X)$;
\item  they  have the same
semi-positive Chern-Ricci curvature, i.e.
$$\Theta^{(1)}(\omega_1)=\Theta^{(1)}(\omega_2)=\Theta^{(1)}(\omega_3)\geq 0;$$
\item they have constant positive  Chern
scalar curvatures. \ed  Moreover,
 \bd \item $\omega_1$ has \textbf{\emph{positive}}  constant Riemannian scalar
 curvature;

 \item $\omega_2$ has \textbf{\emph{{zero}}} Riemannian scalar curvature;

\item $\omega_3$ has \textbf{\emph{negative}} constant Riemannian scalar curvature.
\ed \etheorem

\noindent To the best of our knowledge, this is the first example to
show
 that
 Hermitian manifolds with nonnegative first Chern class can admit  Hermitian metrics with strictly negative constant Riemannian scalar
curvature.\vspace{0.2cm}

We observe the following identity on general compact Hermitian
manifolds.

\bproposition\label{D} On a compact Hermitian manifold $(M,\omega)$,
for any $1\leq k\leq n-1$, we have \beq  \int \sq \p\omega\wedge
\bp\omega \cdot
\frac{\omega^{n-3}}{(n-3)!}=\|\p^*\omega\|^2-\|\p\omega\|^2,\eeq and
\beq \int \sq \omega^{n-k-1}\wedge \p\bp\omega^{k}=
(n-3)!k(n-k-1)\left(\|\p\omega\|^2-\|\p^*\omega\|^2\right)\label{333}.\eeq
\eproposition

\noindent The form $$ \sq \omega^{n-k-1}\wedge \p\bp\omega^{k}$$ was
firstly introduced in \cite{FWW10} by Fu-Wang-Wu  to define a
generalized Gauduchon's metric. More precisely, a metric $\omega$
satisfying $\sq \omega^{n-k-1}\wedge \p\bp\omega^{k}=0$ for $1\leq
k\leq n-1$ is called a $k$-Gauduchon metric. The $(n-1)$-Gauduchon
metric is the original Gauduchon metric. It is well-known that, the
Hopf manifold $\S^{2n+1}\times S^1$ can not support a metric with
$\p\bp\omega=0$ (SKT) or $d^*\omega=0$ (balanced metric). However,
they showed in \cite{FWW10} that on $\S^5\times \S^1$, there exists
a $1$-Gauduchon metric $\omega$, i.e. $\omega\wedge \p\bp\omega=0$.

  A straightforward application of Proposition \ref{D} is the following
interesting fact: \bcorollary\label{K} If $(M,\omega)$ is
$k$-Gauduchon for $1\leq k\leq n-2$ and also balanced, then
$(M,\omega)$ is K\"ahler. \ecorollary \noindent One can also get the
following analogue  in the ``conformal" setting:
\bcorollary\label{Z} On a compact complex manifold, the following
are equivalent:

\bd \item $(M,\omega)$ is conformally K\"ahler;

\item $(M,\omega)$ is conformally $k$-Gauduchon for $1\leq k\leq
n-2$, and conformally balanced.

\ed

In particular, the following are equivalent: \bd \item[(3)]
$(M,\omega)$ is K\"ahler;

\item[(4)] $(M,\omega)$ is  $k$-Gauduchon for $1\leq k\leq
n-2$, and conformally balanced;

\item[(5)]  $(M,\omega)$ is conformally balanced and $\Lambda^2(\p\bp\omega)=0$.
\ed \ecorollary

\noindent We need to point out that the equivalence of $(3)$ and
$(4)$ is also proved by Ivanov and Papadopoulos in
\cite[Theorem~1.3]{IP13} (see also \cite[Proposition~2.4]{FU13}).
 However, in \cite{FWW13}, Fu-Wang-Wu proved that there
exists a non-K\"ahler $3$-fold which can support a $1$-Gauduchon
metric and a balanced metric simultaneously. By Corollary \ref{K},
they must be different Hermitian metrics.
  For more works on Hermitian manifolds with specials metrics,   we refer the reader to \cite{AI01, AT13, AD03, FG, FPS, F1, F2, F3, FY, FLY, FWW10, FWW13, FZ14, GGP, GI01, L, LY12, LY14, Yang2, Yang14, Mi83, T,
  Po1,Po2,Po3} and  the references therein.\vspace{0.2cm}

The paper is organized as follows: In Section \ref{basic}, we
introduce several basic terminologies which will be used frequently
in the paper. In Section \ref{geometry}, we study the  geometry of
the Levi-Civita Ricci curvature and prove Theorem \ref{main-1},
Corollary \ref{c1}, Theorem \ref{main-100} and Corollary \ref{c2}.
In Section \ref{riccirelations}, we investigate a variety of Ricci
curvature and scalar curvature relations over Hermitian manifolds
and establish Theorem \ref{riccire} and Corollary \ref{sphere}. In
Section \ref{special}, we study some special metrics on Hermitian
manifolds and prove Proposition \ref{D}, Corollary \ref{K} and
Corollary \ref{Z}. In Section \ref{example}, we construct various
precise Hermitian metrics on Hermitian manifolds and prove Theorem
\ref{example1} and Theorem \ref{example2}. In Section
\ref{appendix}, we include some straightforward computations on
Hermitian manifolds for the reader's convenience.\\

\textbf{Acknowledgement.} This work was partially supported by
China's Recruitment
 Program of Global Experts and  National Center for Mathematics and Interdisciplinary Sciences,
 Chinese Academy of Sciences. The authors would like to thank J.-X. Fu,
Y.-S. Poon and S. -T. Yau for many useful discussions. The second
author wishes to thank Valentino Tosatti and Ben Weikove for their
invaluable suggestions and support.

\section{Background materials}\label{basic}

\subsection{Ricci curvature on almost Hermitian manifolds}
Let $(M^{}, g, \nabla)$ be a $2n$-dimensional Riemannian manifold
with Levi-Civita connection $\nabla$. The tangent bundle of $M$ is
denoted by $T_\R M$. The curvature tensor of $(M,g,\nabla)$  is
defined as \beq
R(X,Y,Z,W)=g\left(\nabla_X\nabla_YZ-\nabla_Y\nabla_XZ-\nabla_{[X,Y]}Z,W\right)\eeq
for any $X,Y,Z,W\in T_\R M$. Let $T_\C M=T_\R M\ts \C$  be the
complexification of the tangent bundle $T_\R M$. We can extend the
metric $g$,  the Levi-Civita connection $\nabla$ to $T_{\C}M$ in the
$\C$-linear way. For instance, for any $a,b\in \C$ and $X,Y\in T_\C
M$, \beq g(aX, bY):=ab\cdot g(X,Y).\eeq Hence for any $a,b,c,d\in
\C$ and $X,Y,Z,W\in T_\C M$,  \beq R(aX,bY,cZ, dW)=abcd\cdot
R(X,Y,Z,W).\eeq

\noindent Let $(M,g,J)$ be an almost Hermitian manifold, i.e.,
$J:T_\R M\>T_\R M$ with $J^2=-1$, and
 for any $X,Y\in T_\R M$, $ g(JX,JY)=g(X,Y)$. We can also extend $J$ to  $T_\C M$ in the
$\C$-linear way. Hence for any $X,Y\in T_\C M$, we still have \beq
g(JX,JY)=g(X,Y).\eeq

Let $\{x^1,\cdots, x^n,x^{n+1},\cdots， x^{2n}\}$ be the local real
coordinates on the almost Hermitian manifold $(M,J,g)$. In order to
use Einstein summations, we use the following convention: \beq
\{x^i\} \qtq{for} 1\leq i\leq n;\ \ \ \ \  \{x^I\} \qtq{for} n+1\leq
I\leq 2n \qtq{and}\ \ \ I=i+n.\eeq Moreover, we assume, \beq
J\left(\frac{\p}{\p x^i}\right)=\frac{\p}{\p x^I} \qtq{and}
J\left(\frac{\p}{\p x^I}\right)=-\frac{\p}{\p x^i}.\eeq By using
real coordinates $\{x^i,x^I\}$, the Riemannian metric is represented
by
$$ds_g^2=g_{i\ell}dx^i\ts dx^\ell+g_{iL}dx^i\ts dx^L+g_{I\ell} dx^I\ts dx^\ell+g_{IJ}dx^I\ts dx^J,$$
where the metric components $g_{i\ell}, g_{iL}, g_{I\ell}$ and
$g_{IL}$ are defined in the obvious way by using $\frac{\p}{\p x^i},
\frac{\p}{\p x^\ell}, \frac{\p}{\p x^I}, \frac{\p}{\p x^L}$. By the
$J$-invariant property of the metric $g$, we have \beq
g_{i\ell}=g_{IL}, \qtq{and} g_{iL}=g_{Li}=-g_{\ell
I}=-g_{I\ell}.\label{1001} \eeq

\noindent We also use \emph{complex coordinates} $\{z^i,\bar
z^i\}_{i=1}^n$ on $M$: \beq z^i:=x^i+\sq x^I,\ \ \ \bar z^i:=x^i-\sq
x^I.\eeq (Note that if the almost complex structure $J$ is
integrable, $\{z^i\}_{i=1}^n$ are the local holomorphic
coordinates.) We  define, for $1\leq i\leq n$, \beq dz^i:=dx^i+\sq
dx^I,\ \ \ d\bar z^i:=dx^i-\sq dx^I\eeq and\beq \frac{\p}{\p
z^i}:=\frac{1}{2}\left(\frac{\p}{\p x^i}-\sq \frac{\p}{\p
x^I}\right),\ \ \ \frac{\p}{\p\bar
z^i}:=\frac{1}{2}\left(\frac{\p}{\p x^i}+\sq \frac{\p}{\p
x^I}\right).\eeq Therefore,  \beq \frac{\p}{\p x^i}=\frac{\p}{\p
z^i}+\frac{\p}{\p \bar z^i},\ \ \ \frac{\p}{\p
x^I}=\sq\left(\frac{\p}{\p z^i}-\frac{\p}{\p \bar
z^i}\right).\label{1002}\eeq By the $\C$-linear extension,\beq
J\left(\frac{\p}{\p z^i}\right)=\sq \frac{\p}{\p z^i},\ \ \
J\left(\frac{\p}{\p \bar z^i}\right)=-\sq \frac{\p}{\p \bar
z^i}.\eeq Let's define a Hermitian form $h:T_\C M\times T_\C M\>\C$
by \beq h(X,Y):= g(X,Y),
 \ \ \ \ \ X,Y\in T_\C M.\label{complex}\eeq
By $J$-invariant property of $g$, \beq h_{ij}:=h\left(\frac{\p}{\p
z^i},\frac{\p}{\p z^j}\right)=0, \qtq{and} h_{\bar i\bar
j}:=h\left(\frac{\p}{\p \bar z^i},\frac{\p}{\p \bar
z^j}\right)=0\eeq and \beq h_{i\bar j}:=h\left(\frac{\p}{\p
z^i},\frac{\p}{\p \bar z^j}\right)=\frac{1}{2}\left(g_{ij}+\sq
g_{iJ}\right).\label{1003}\eeq It is obvious that $(h_{i\bar j})$ is
a positive Hermitian  matrix. Here, we always use the convention
$h_{\bar j i}=h_{i\bar j}$ since $g$ is symmetric over $T_\C M$. The
(transposed) inverse matrix of $(h_{i\bar j})$ is denoted by
$(h^{i\bar j})$, i.e. $h^{i\bar \ell}\cdot h_{k\bar
\ell}=\delta_{k}^i$. One can also show that, \beq h^{i\bar
j}=2\left(g^{ij}-\sq g^{iJ}\right)\label{1004}\eeq

\noindent From the definition equation (\ref{complex}), we see the
following well-known relation on $T_\C M$: \beq
ds_h^2=\frac{1}{2}ds_g^2-\frac{\sq}{2} \omega\eeq where $\omega$ is
the fundamental $2$-form associated to the $J$-invariant metric $g$:
 \beq \omega(X,Y)=g(JX,Y).\eeq
   In local complex coordinates,
 \beq \omega=\sq h_{i\bar j} dz^i\wedge d\bar z^j.\eeq

 \noindent
In the following, we shall use the components of the complexified
curvature tensor $R$, for example,  \beq R_{i\bar j k\bar
\ell}:=R\left(\frac{\p}{\p z^i}, \frac{\p}{\p \bar z^j},
\frac{\p}{\p z^k}, \frac{\p}{\p \bar z^\ell}\right),\ \ \ \eeq and
in particular we use the following notation for the complexified
curvature tensor: \beq R_{i j k \ell}:=R\left(\frac{\p}{\p z^i},
\frac{\p}{\p z^j}, \frac{\p}{\p z^k}, \frac{\p}{\p z^\ell}\right).\
\ \ \eeq It is obvious that, the components of ($\C$-linear)
complexified curvature tensor have the same properties as the
components of the real curvature tensor. We list some properties of
$R_{i\bar j k\bar \ell}$ for examples:\beq R_{i\bar j
k\ell}=-R_{\bar j i k\ell},\ \ \ R_{i\bar j k\bar\ell}=R_{k\bar \ell
i\bar j},\eeq and in particular, the (first) Bianchi identity holds:
\beq R_{i\bar j k\bar \ell}+R_{ik\bar\ell \bar j}+R_{i\bar \ell \bar
jk}=0.\label{bianchi}\eeq

\bdefinition  Let $\{e_i\}_{i=1}^{2n}$ be a local orthonormal basis
of $(T_\R M,g)$, the Riemannian Ricci curvature of $(M,g)$ is \beq
Ric (X,Y):=\sum_{i=1}^{2n}R(e_i,X,Y,e_i), \eeq and the corresponding
Riemannian scalar curvature is \beq
s=\sum_{j=1}^{2n}Ric(e_i,e_i).\eeq

\edefinition

\blemma\label{key} On an almost Hermitian manifold $(M,h)$,  the
Riemannian Ricci curvature of the Riemannian manifold $(M,g)$
satisfies \beq Ric(X,Y)=h^{i\bar \ell}\left[R\left(\frac{\p}{\p
z^i}, X,Y, \frac{\p}{\p\bar z^\ell}\right)+R\left(\frac{\p}{\p z^i},
Y,X, \frac{\p}{\p\bar z^\ell}\right)\right]\label{ricci}\eeq for any
$X,Y\in T_\R M$. The Riemannian scalar curvature is \beq s=2h^{i\bar
j}h^{k\bar \ell}\left(2R_{i\bar \ell k\bar j}-R_{i\bar j k\bar
\ell}\right).\label{scalar}\eeq \elemma \bproof See Lemma \ref{Akey}
in the Appendix. \eproof

\noindent In order to formulate the curvature relations more
effectively, we introduce new curvature notations as following:

\bdefinition The Riemannian Ricci tensor can also be extended to
$T_\C M$, and we denote the associated $(1,1)$-form component by
\beq \mathscr Ric=\sq \mathscr R_{i\bar j}dz^i\wedge d\bar z^j
\qtq{with} \mathscr R_{i\bar j}:=Ric\left(\frac{\p}{\p z^i},
\frac{\p}{\p\bar z^j}\right). \eeq
 The \emph{Hermitian-Ricci curvature} of the
complexified curvature tensor is \beq Ric_H=\sq R_{i\bar
j}dz^i\wedge d\bar z^j, \qtq{with} R_{i\bar j}:=h^{k\bar
\ell}R_{i\bar j k\bar \ell}.\eeq
 The  corresponding \emph{Hermitian scalar curvature} of $h$ is given by\beq
s_H=h^{i\bar j}R_{i\bar j}. \eeq We also define the \emph{Riemannian
type scalar curvature} as \beq s_R=h^{i\bar \ell}h^{k\bar j}R_{i\bar
j k\bar \ell}.\eeq It is obvious that $s_H\neq s_R$ in general.
Similarly, we can define the third and fourth Hermitian-Ricci
curvatures $R^{(3)}$ and $R^{(4)}$ respectively,
$$ R^{(3)}=\sq  R^{(3)}_{i\bar j} dz^i\wedge d\bar z^j\qtq{with}  R^{(3)}_{i\bar j}=h^{k\bar\ell}  R_{i\bar \ell k \bar j},$$
and
$$ R^{(4)}=\sq  R^{(4)}_{i\bar j} dz^i\wedge d\bar z^j\qtq{with}  R^{(4)}_{i\bar j}=h^{k\bar\ell}  R_{k\bar j i\bar \ell }.$$
It is easy to see that $R^{(3)}=R^{(4)}$.

\edefinition

\bcorollary On an almost Hermitian manifold $(M,h)$, we have \beq
\mathscr{R}_{i\bar j}=2\left(h^{k\bar \ell}R_{k \bar j i\bar
\ell}\right) -R_{i\bar j}\label{ricci2} \eeq and \beq s=2h^{i\bar
j}\mathscr R_{i\bar j}=4s_R-2s_H.\label{scalar2}\eeq \ecorollary
\bproof  (\ref{ricci2}) follows from (\ref{ricci}) and the Bianchi
identity (\ref{bianchi}). (\ref{scalar2}) follows from
(\ref{scalar}). \eproof

\subsection{Curvatures on Hermitian manifolds}
Let $(M, h, J)$ be an almost Hermitian manifold. The Nijenhuis
tensor $N_J:\Gamma(M,T_\R M)\times \Gamma(M,T_\R M)\>\Gamma(M,T_\R
M)$ is defined as \beq N_J(X,Y)=[X,Y]+J[JX,Y]+J[X,JY]-[JX,JY].\eeq
The almost complex structure $J$ is called \emph{integrable} if
$N_J\equiv 0$ and then we call $(M,g,J)$ a Hermitian manifold. By
Newlander-Nirenberg's theorem, there exists a real coordinate system
$\{x^i,x^I\}$ such that $z^i=x^i+\sq x^I$ are  local holomorphic
coordinates on $M$. Moreover, we have $T_\C M=T^{1,0}M\ds T^{0,1}M$
where $$T^{1,0}M=span_\C\left\{\frac{\p}{\p z^1},\cdots,\frac{\p}{\p
z^n}\right\}\qtq{and}T^{0,1}M=span_\C\left\{\frac{\p}{\p \bar
z^1},\cdots,\frac{\p}{\p \bar z^n}\right\}.$$

\noindent Let $\phi$ be a $(p,q)$-form on $(M,g)$, and \beq
\phi=\frac{1}{p!q!}\sum_{i_1,\cdots, i_p, j_1,\cdots,
j_q}\phi_{i_1\cdots i_p \bar{j_1}\cdots \bar{j_q}}dz^{i_1}\wedge
\cdots\wedge  dz^{i_p}\wedge d\bar z^{j_1}\wedge \cdots\wedge d\bar
z^{j_q},\eeq where $\phi_{i_1\cdots i_p \bar{j_1}\cdots \bar{j_q}}$
is skew symmetric in $i_1,\cdots, i_p$ and also skew symmetric in
$j_1,\cdots, j_q$. The local inner product is defined as \beq
|\phi|^2= \la\phi,\phi\ra=\frac{1}{p!q!}h^{i_1\bar {\ell_1}}\cdots
h^{i_p\bar {\ell_p}}h^{k_1\bar {j_1}}\cdots h^{k_q\bar {j_q}}
\phi_{i_1\cdots i_p \bar{j_1}\cdots \bar{j_q}}\cdot
\bar{\phi_{\ell_1\cdots \ell_p \bar{k_1}\cdots \bar{k_q}}}.\eeq The
norm on $\Om^{p,q}(M)$ is \beq \|\phi\|^2=(\phi,\phi)=\int
\la\phi,\phi\ra\frac{\omega^n}{n!}.\eeq It is well-known that there
exists a real isometry $*:\Om^{p,q}(M)\>\Om^{n-q,n-p}(M)$ such that
\beq (\phi,\psi)=\int \phi\wedge *\bar\psi,\eeq for $\phi,\psi\in
\Om^{p,q}(M)$.

\section{Geometry of the Levi-Civita Ricci
curvature}\label{geometry}

\subsection{The  Levi-Civita connection and Chern connection on $(T^{1,0}M,h)$}

\subsubsection{The induced Levi-Civita connection on $(T^{1,0}M,h)$}
\noindent Since $T^{1,0}M$ is a subbundle of $T_{\C}M$, there is an
induced connection $\hat\nabla$ on $T^{1,0}M$ given by \beq
\hat\nabla=\pi\circ\nabla:
\Gamma(M,T^{1,0}M)\stackrel{\nabla}{\rightarrow}\Gamma(M, T_{\C}M\ts
T_{\C}M)\stackrel{\pi}{\rightarrow}\Gamma(M,T_{\C}M\ts T^{1,0}M).
\label{metricconnection}\eeq Moreover, $\hat \nabla$ is  a metric
connection on the Hermitian holomorphic vector bundle $(T^{1,0}M,
h)$ and it is determined by the  relations \beq
\hat\nabla_{\frac{\p}{\p z^i}}\frac{\p}{\p
z^k}:=\Gamma_{ik}^p\frac{\p}{\p z^p} \qtq{and}
\hat\nabla_{\frac{\p}{\p \bar z^j}}\frac{\p}{\p z^k}:=\Gamma_{\bar
jk}^p\frac{\p}{\p z^p} \eeq where \beq
\Gamma_{ij}^k=\frac{1}{2}h^{k\bar \ell}\left(\frac{\p h_{j\bar
\ell}}{\p z^i}+\frac{\p h_{i\bar \ell}}{\p z^j}\right), \qtq{and}
\Gamma_{\bar i j}^k=\frac{1}{2} h^{k\bar \ell}\left( \frac{\p
h_{j\bar\ell}}{\p\bar z^i}-\frac{\p h_{j\bar i}}{\p\bar
z^\ell}\right).\eeq

\noindent The curvature  tensor $\mathfrak{R}\in \Gamma(M,\Lambda^2
T_{\C}M\ts T^{*1,0}M\ts T^{1,0}M)$ of $\hat\nabla$ is given by \beq
\mathfrak{R}(X,Y)s
=\hat\nabla_{X}\hat\nabla_Ys-\hat\nabla_Y\hat\nabla_Xs-\hat\nabla_{[X,Y]}s\eeq
for any $X,Y\in T_{\C}M$ and $s\in T^{1,0}M$. The curvature tensor
$\mathfrak{R}$ has components \beq \mathfrak{R}_{i\bar
jk}^{\ell}=-\left(\frac{\p \Gamma^{\ell}_{ik}}{\p \bar z^j}-\frac{\p
\Gamma^{\ell}_{\bar jk}}{\p z^i}+\Gamma_{ ik}^{s}\Gamma^{\ell}_{\bar
js}-\Gamma_{ \bar jk}^{s}\Gamma^{\ell}_{s i}\right)\qtq{and}
 \mathfrak{R}_{\bar j ik}^{\ell}=-\mathfrak{R}_{i\bar j k}^\ell;\label{curv}\eeq
 \beq \mathfrak{R}_{i jk}^{\ell}=-\left(\frac{\p
\Gamma^{\ell}_{ ik}}{\p z^j}-\frac{\p \Gamma^{\ell}_{jk}}{\p
z^i}+\Gamma_{ik}^{s}\Gamma^{\ell}_{s j}-\Gamma_{
jk}^{s}\Gamma^{\ell}_{s i}\right);\eeq and \beq \mathfrak{R}_{\bar
i\bar j k}^{\ell}=-\left(\frac{\p \Gamma^{\ell}_{\bar ik}}{\p \bar
z^j}-\frac{\p \Gamma^{\ell}_{\bar jk}}{\p \bar z^i}+\Gamma_{\bar
ik}^{s}\Gamma^{\ell}_{\bar js}-\Gamma_{\bar
jk}^{s}\Gamma^{\ell}_{\bar is}\right).\eeq

\noindent With respect to the Hermitian metric $h$ on $T^{1,0}M$, we
use the  convention  \beq \mathfrak{R}_{\bullet\bullet k\bar
\ell}:=\sum_{s=1}^n\mathfrak{R}_{\bullet\bullet k}^{s}h_{s\bar
\ell}. \eeq

\bcorollary[{\cite[Proposition~2.1]{LY12}}] We have the following
relations: $$R_{ijk}^\ell=\mathfrak{R}_{ijk}^\ell,\ \ \ R_{\bar
i\bar j k}^{\ell}=\mathfrak{R}_{\bar i\bar j k}^\ell,$$ and \beq
R_{i\bar j k}^\ell=-\left(\frac{\p \Gamma^{\ell}_{ik}}{\p \bar
z^j}-\frac{\p \Gamma^{\ell}_{\bar jk}}{\p z^i}+\Gamma_{
ik}^{s}\Gamma^{\ell}_{\bar js}-\Gamma_{ \bar jk}^{s}\Gamma^{\ell}_{s
i}-\Gamma_{\bar s i}^\ell\cdot\bar{\Gamma_{\bar k
j}^s}\right)=\mathfrak{R}_{i\bar j k}^\ell+\Gamma_{\bar s
i}^\ell\cdot\bar{\Gamma_{\bar k j}^s}.\label{key0}\eeq \ecorollary

\noindent Next, we define Ricci curvatures and scalar curvatures for
$(T^{1,0}M, h, \hat\nabla)$.

 \bdefinition The \emph{first Levi-Civita Ricci curvature}  of
the Hermitian  vector bundle $\left(T^{1,0}M, h, \hat \nabla\right)$
is
 \beq \mathfrak{R}^{(1)}=\sq\mathfrak
R^{(1)}_{i\bar j}dz^i\wedge d\bar z^j \qtq{with}
\mathfrak{R}^{(1)}_{i\bar j}=h^{k\bar \ell} \mathfrak{R}_{i\bar j
k\bar\ell} \eeq  and the \emph{second Levi-Civita Ricci curvature}
of it is
 \beq \mathfrak{R}^{(2)}=\sq \mathfrak{R}^{(2)}_{i\bar j}dz^i\wedge d\bar z^j \qtq{with} \mathfrak{R}^{(2)}_{i\bar j}=h^{k\bar \ell} \mathfrak{R}_{ k\bar
 \ell i\bar j}.\label{2riccilc}
\eeq The \emph{Levi-Civita scalar curvature} of $\hat\nabla$ on
$T^{1,0}M$ is denoted by \beq s_{LC}=h^{i\bar j}h^{k\bar
\ell}\mathfrak{R}_{i\bar j k\bar \ell}. \eeq Similarly, we can
define $\mathfrak R^{(3)}$ and $\mathfrak R^{(4)}$ as
$$\mathfrak R^{(3)}=\sq \mathfrak R^{(3)}_{i\bar j} dz^i\wedge d\bar z^j\qtq{with} \mathfrak R^{(3)}_{i\bar j}=h^{k\bar\ell} \mathfrak R_{i\bar \ell k \bar j},$$
and
$$\mathfrak R^{(4)}=\sq \mathfrak R^{(4)}_{i\bar j} dz^i\wedge d\bar z^j\qtq{with} \mathfrak R^{(4)}_{i\bar j}=h^{k\bar\ell} \mathfrak R_{k\bar j i\bar \ell }.$$
 \edefinition
\vskip 1\baselineskip
\subsubsection{Curvature of the Chern connection on $(T^{1,0}M,h)$} On the
Hermitian holomorphic vector bundle $(T^{1,0}M,h)$, the Chern
connection $\nabla^{Ch}$ is the unique connection which is
compatible with the complex structure and also the Hermitian metric.
The curvature tensor of $\nabla^{Ch}$ is denoted by $\Theta$ and its
curvature components are \beq \Theta_{i\bar j k\bar
\ell}=-\frac{\p^2 h_{k\bar \ell}}{\p z^i\p \bar z^j}+h^{p\bar
q}\frac{\p h_{p\bar \ell}}{\p \bar z^j}\frac{\p h_{k\bar q}}{\p
z^i}.\label{cherncurvature} \eeq
 It is
well-known that the \emph{(first) Chern Ricci curvature} \beq
\Theta^{(1)}:= \sq\Theta^{(1)}_{i\bar j} dz^i\wedge d\bar z^j\eeq
represents the first Chern class $ c_1(M)$ of $M$ where
  \beq
\Theta^{(1)}_{i\bar j}= h^{k\bar \ell}\Theta_{i\bar j k\bar \ell}
=-\frac{\p^2 \log \det(h_{k\bar \ell})}{\p z^i\p\bar
z^j}.\label{1chern}\eeq The \emph{second Chern Ricci curvature}
$\Theta^{(2)}=\sq \Theta^{(2)}_{i\bar j}dz^i\wedge d\bar z^j$ with
components \beq \Theta^{(2)}_{i\bar j}=h^{k\bar \ell}\Theta_{k\bar
\ell i\bar j}. \label{2chern}\eeq The \emph{Chern scalar curvature}
$s_C$ of the Chern curvature tensor $\Theta$ is defined by \beq
s_{C}=h^{i\bar j}h^{k\bar \ell}\Theta_{i\bar j k\bar \ell}.\eeq
Similarly, we can define $\Theta^{(3)}$ and $\Theta^{(4)}$.

\vskip 1\baselineskip

 The (first) Chern-Ricci curvature $\Theta^{(1)}$ represents
the first Chern class $ c_1(M)$, but in general, the first
Levi-Civita Ricci curvature $\mathfrak{R}^{(1)}$ does not represent
a class in $H_{dR}^2(M)$ or $H^{1,1}_{\bp}(M)$, since it is not
$d$-closed. We shall
 explore more geometric properties of $\fR^{(1)}$ in the following sections.

\subsection{Elementary computations on Hermitian manifolds} In this
subsection, we recall some elementary and well-known computational
lemmas on Hermitian manifolds.

\blemma Let $(M,h)$ be a compact Hermitian manifold and $\omega=\sq
h_{i\bar j}dz^i\wedge d\bar z^j$.
 \beq\p^*\omega=-\sq
\Lambda\left(\bp\omega\right)=-2\sq \Gamma_{\bar j k}^kd\bar z^j
\qtq{and} \bp^*\omega=\sq \Lambda\left(\p\omega\right)=2\sq
\bar{\Gamma_{\bar ik}^k}dz^i.\label{key1}\eeq \elemma \bproof By the
well-known Bochner formula (e.g. \cite{LY12}), $$[\bp^*,L]=\sq
\left(\p+\tau\right)$$ where $\tau=[\Lambda,\p\omega]$, we see
$\bp^*\omega=\sq \Lambda\left(\p\omega\right)=2\sq \bar{\Gamma_{\bar
ik}^k}dz^i$. \eproof

 \blemma Let $(M,h,\omega)$ be a Hermitian
manifold. For any $p\in M$, there exist local holomorphic ``normal
coordinates" $\{z^i\}$ centered at $p$ such that \beq h_{i\bar
j}(p)=\delta_{ij} \qtq{and} \Gamma_{ij}^k(p)=0.\label{coordinate}
\eeq In particular, at $p$, we have \beq \Gamma_{\bar j
i}^k=\frac{\p h_{i\bar k}}{\p \bar z^j}=-\frac{\p h_{i\bar
j}}{\p\bar z^k}.\eeq
 \elemma

Let $T$ be the torsion tensor of the Hermitian metric $\omega$, i.e.
\beq T_{ij}^k=h^{k\bar\ell}\left(\frac{\p h_{j\bar \ell}}{\p
z^i}-\frac{\p h_{i\bar \ell}}{\p z^j}\right).\eeq In the following,
we shall use the  conventions: \beq T\boxdot \bar T:=h^{p\bar q}
h_{k \bar \ell}T_{ip}^k\cdot \bar {T_{jq}^\ell} dz^i\wedge d\bar
z^j, \qtq{and} T\circ \bar T:=h^{p\bar q}h^{s\bar t} h_{k \bar j}
h_{i \bar \ell}T_{sp}^k\cdot \bar {T_{tq}^\ell} dz^i\wedge d\bar
z^j.\eeq

\noindent It is obvious that the $(1,1)$-forms $T\circ \bar T$ and
$T\boxdot \bar T$ are not the same.

\blemma At a fixed point $p$ with ``normal coordinates"
(\ref{coordinate}), we have \beq (T\boxdot\bar T)_{i\bar
j}=T_{ip}^k\cdot \bar{T_{jp}^k}=4\sum_{p,k}\frac{\p h_{p\bar k}}{\p
z^i}\cdot {\frac{\p h_{ k\bar p}} {\p \bar z^j}}\eeq and \beq
(T\circ\bar T)_{i\bar j}=T_{pq}^j\cdot
\bar{T_{pq}^i}=4\sum_{p,q}\frac{\p h_{q\bar j}}{\p z^p}\cdot
{\frac{\p h_{ i\bar q}} {\p \bar z^p}}.\eeq Moreover
$$tr_\omega \left(\sq T\boxdot\bar T\right)=tr_\omega\left(\sq T\circ \bar T\right)=|T|^2.$$
\elemma

\blemma At a fixed point $p$ with ``normal coordinates"
(\ref{coordinate}), we have the $(1,1)$ form\beq
T((\p^*\omega)^{\#})=-4\sq \frac{\p h_{i\bar j}}{\p z^s}\frac{\p
h_{\ell\bar \ell}}{\p \bar z^s} dz^i\wedge d\bar z^j,\label{98}\eeq
where $(\p^*\omega)^{\#}$ is the dual vector of  the $(0,1)$-form
$\p^*\omega$. \elemma \bproof Since $\p^*\omega=-2\sq \Gamma_{\bar s
\ell}^\ell d\bar z^s$, the corresponding $(1,0)$ type vector field
is \beq (\p^*\omega)^\#=-2{\sq}h^{i\bar s}\Gamma_{\bar s \ell}^\ell
\frac{\p}{\p z^i}.\eeq The $(1,1)$ form \beq
T((\p^*\omega)^{\#})=-2\sq h_{k\bar j}T^k_{pi}\left(h^{p\bar
s}\Gamma_{\bar s \ell}^\ell \right)dz^i\wedge d\bar
z^j=-4\sq\sum_{\ell, s, i, j} \frac{\p h_{i\bar j}}{\p z^s}\frac{\p
h_{\ell\bar \ell}}{\p \bar z^s} dz^i\wedge d\bar z^j.\eeq \eproof

\blemma At a fixed point $p$ with ``normal coordinates"
(\ref{coordinate}), we have
$$\p\p^*\omega=\sq(\p\p^*\omega)_{i\bar j}dz^i\wedge d\bar z^j$$ where \beq
(\p\p^*\omega)_{i\bar j}=\sum_q\left(\frac{\p^2 h_{q\bar j}}{\p
z^i\p\bar z^q}-\frac{\p^2 h_{q\bar q}}{\p z^i\p\bar z^j}\right)
+\frac{1}{2}(T\boxdot\bar T)_{i\bar j}.\label{dd*o}\eeq Similarly,
$$\bp\bp^*\omega=\sq(\bp\bp^*\omega)_{i\bar j}dz^i\wedge d\bar z^j$$ where \beq
(\bp\bp^*\omega)_{i\bar j}=\sum_q\left(\frac{\p^2 h_{i\bar q}}{\p
z^q\p\bar z^j}-\frac{\p^2 h_{q\bar q}}{\p z^i\p\bar z^j}\right)
+\frac{1}{2}(T\boxdot\bar T)_{i\bar j}.\label{pp*o}\eeq \elemma

\bproof Since $$\p^*\omega=-2\sq \Gamma_{\bar s \ell}^\ell d\bar
z^s=-\sq h^{\ell \bar q}\left(\frac{\p h_{\ell\bar q}}{\p\bar
z^s}-\frac{\p h_{\ell\bar s}}{\p\bar z^q}\right)d\bar z^s,$$ we have
\be \p\p^*\omega&=&-\sq h^{\ell\bar q}\left(\frac{\p^2 h_{\ell\bar
q}}{\p z^i\p\bar z^s}-\frac{\p^2 h_{\ell\bar s}}{\p z^i\p\bar
z^q}\right)dz^i\wedge d\bar z^s-\sq \frac{\p h^{\ell \bar q}}{\p
z^i}\left(\frac{\p h_{\ell\bar q}}{\p\bar z^s}-\frac{\p h_{\ell\bar
s}}{\p\bar z^q}\right)dz^i\wedge d\bar z^s\\
&=& \sq\sum_q\left(\frac{\p^2 h_{q\bar s}}{\p z^i\p\bar
z^q}-\frac{\p^2 h_{q\bar q}}{\p z^i\p\bar z^s}\right)dz^i\wedge
d\bar z^s+2\sq \sum_{q,\ell}\frac{\p h_{q\bar \ell}}{\p z^i}\frac{\p
h_{\ell\bar q}}{\p\bar z^s} dz^i\wedge d\bar z^s\\
&=& \sq\sum_q\left(\frac{\p^2 h_{q\bar s}}{\p z^i\p\bar
z^q}-\frac{\p^2 h_{q\bar q}}{\p z^i\p\bar z^s}\right)dz^i\wedge
d\bar z^s+\frac{\sq}{2}T\boxdot\bar T.\ee \eproof

The following lemmas follow from straightforward computations.
\blemma At a fixed point $p$ with ``normal coordinates"
(\ref{coordinate}), we have \beq \mathfrak{R}_{i\bar j k\bar
\ell}=-\frac{1}{2}\left(\frac{\p^2 h_{i\bar \ell}}{\p z^k\p\bar
z^j}+\frac{\p^2 h_{k\bar j}}{\p z^i\p\bar
z^\ell}\right)-\sum_q\frac{\p h_{q\bar \ell}}{\p z^i}\frac{\p
h_{k\bar q}}{\p\bar z^j}, \eeq \beq R_{i\bar j k\bar
\ell}=\mathfrak{R}_{i\bar j k\bar \ell}-\sum_q\frac{\p h_{q\bar
j}}{\p z^k}\frac{\p h_{i\bar q}}{\p\bar
z^\ell}=-\frac{1}{2}\left(\frac{\p^2 h_{i\bar \ell}}{\p z^k\p\bar
z^j}+\frac{\p^2 h_{k\bar j}}{\p z^i\p\bar
z^\ell}\right)-\sum_q\left(\frac{\p h_{q\bar \ell}}{\p z^i}\frac{\p
h_{k\bar q}}{\p\bar z^j}+\frac{\p h_{q\bar j}}{\p z^k}\frac{\p
h_{i\bar q}}{\p\bar z^\ell}\right) \eeq and \beq \Theta_{i\bar j
k\bar\ell}=-\frac{\p^2 h_{k\bar \ell}}{\p z^i\p \bar
z^j}+\sum_q\frac{\p h_{q\bar \ell}}{\p \bar z^j}\frac{\p h_{k\bar
q}}{\p z^i}.\eeq

\elemma

\blemma At a fixed point $p$ with ``normal coordinates"
(\ref{coordinate}), we have \beq \mathfrak{R}^{(1)}_{i\bar j}=\sum_k
\mathfrak{R}_{i\bar j k \bar k}=-\frac{1}{2}\sum_k\left(\frac{\p^2
h_{i\bar k}}{\p z^k\p\bar z^j}+\frac{\p^2 h_{k\bar j}}{\p z^i\p\bar
z^k}\right)-\frac{1}{4}(T\boxdot \bar T)_{i\bar j};
\label{mathfrakR1}\eeq \beq \mathfrak{R}^{(2)}_{i\bar j}=\sum_k
\mathfrak{R}_{k \bar ki\bar j }=-\frac{1}{2}\sum_k\left(\frac{\p^2
h_{i\bar k}}{\p z^k\p\bar z^j}+\frac{\p^2 h_{k\bar j}}{\p z^i\p\bar
z^k}\right)-\frac{1}{4}(T\circ T)_{i\bar j}; \label{mathfrakR2}\eeq
\beq \mathfrak{R}^{(3)}_{i\bar j}=\sum_k \mathfrak{R}_{i\bar k k
\bar j}=-\frac{1}{2}\sum_k\left(\frac{\p^2 h_{k\bar k}}{\p z^i\p\bar
z^j}+\frac{\p^2 h_{i\bar j}}{\p z^k\p\bar
z^k}\right)+\frac{1}{4}(T(\p^*\omega)^\#)_{i\bar
j};\label{mathfrakR3}\eeq \beq \mathfrak{R}^{(4)}_{i\bar j}=\sum_k
\mathfrak{R}_{k \bar ji\bar k }=-\frac{1}{2}\sum_k\left(\frac{\p^2
h_{k\bar k}}{\p z^i\p\bar z^j}+\frac{\p^2 h_{i\bar j}}{\p z^k\p\bar
z^k}\right)+\frac{1}{4}\bar{(T(\p^*\omega)^\#)}_{i\bar
j}.\label{mathfrakR4}\eeq

\elemma

\blemma At a fixed point $p$ with ``normal coordinates"
(\ref{coordinate}), we have \begin{eqnarray} R^{(1)}_{k\bar
\ell}=R^{(2)}_{k\bar\ell}\nonumber&=&-\frac{1}{2}\sum_s\left(\frac{\p^2
h_{s\bar \ell}}{\p z^k\p\bar z^s}+\frac{\p^2 h_{k\bar s}}{\p
z^s\p\bar z^\ell}\right)-\sum_{q,s}\left(\frac{\p h_{q\bar \ell}}{\p
z^s}\frac{\p h_{k\bar q}}{\p\bar z^s}+\frac{\p h_{k\bar q}}{\p
z^s}\frac{\p h_{q\bar \ell}}{\p\bar z^s}\right)\\
&=&-\frac{1}{2}\sum_s\left(\frac{\p^2 h_{s\bar \ell}}{\p z^k\p\bar
z^s}+\frac{\p^2 h_{k\bar s}}{\p z^s\p\bar
z^\ell}\right)-\frac{(T\circ \bar T)_{k\bar \ell}+(T\boxdot \bar
T)_{k\bar \ell}}{4};\label{R1R2}\end{eqnarray}  \begin{eqnarray}
R^{(3)}_{k\bar\ell}=R^{(4)}_{k\bar\ell}\nonumber&=&-\frac{1}{2}\sum_s\left(\frac{\p^2
h_{k\bar \ell}}{\p z^s\p\bar z^s}+\frac{\p^2 h_{s\bar s}}{\p
z^k\p\bar z^\ell}\right)-\sum_{q,s}\left(\frac{\p h_{q\bar \ell}}{\p
z^k}\frac{\p h_{s\bar q}}{\p \bar z^s}+\frac{\p h_{q\bar s}}{\p
z^s}\frac{\p h_{k\bar q}}{\p \bar z^\ell}\right)\\
&=&-\frac{1}{2}\sum_s\left(\frac{\p^2 h_{k\bar \ell}}{\p z^s\p\bar
z^s}+\frac{\p^2 h_{s\bar s}}{\p z^k\p\bar
z^\ell}\right)+\frac{1}{4}{(T(\p^*\omega)^\#)}_{i\bar
j}+\frac{1}{4}\bar{(T(\p^*\omega)^\#)}_{i\bar j}. \label{R3R4}
\end{eqnarray}

\elemma

\blemma At a fixed point $p$ with ``normal coordinates"
(\ref{coordinate}), we have \beq \Theta^{(1)}_{i\bar
j}=-\sum_q\frac{\p^2 h_{q\bar q}}{\p z^i\p\bar
z^j}+\frac{1}{4}(T\boxdot\bar T)_{i\bar j};\label{theta1}\eeq \beq
\Theta^{(2)}_{i\bar j}=-\sum_q\frac{\p^2 h_{i\bar j}}{\p z^q\p\bar
z^q}+\frac{1}{4}(T\boxdot\bar T)_{i\bar j};\label{theta2}\eeq \beq
\Theta^{(3)}_{i\bar j}=-\sum_q \frac{\p^2 h_{q\bar j}}{\p z^i\p\bar
z^q}-\frac{1}{4}(T\boxdot\bar T)_{i\bar j};\label{theta3}\eeq \beq
\Theta^{(4)}_{i\bar j}=-\sum_q \frac{\p^2 h_{i\bar q}}{\p z^q\p\bar
z^j}-\frac{1}{4}(T\boxdot\bar T)_{i\bar j}.\label{theta4}\eeq

\elemma

\blemma At a fixed point $p$ with ``normal coordinates"
(\ref{coordinate}),  the complexified Ricci curvature is
\begin{eqnarray} \mathscr R_{k\bar
\ell}&=&\nonumber\frac{1}{2}\sum_s\left(\frac{\p^2 h_{s\bar
\ell}}{\p z^k\p\bar z^s}+\frac{\p^2 h_{k\bar s}}{\p z^s\p\bar
z^\ell}\right)- \sum_s\left(\frac{\p^2 h_{k\bar \ell}}{\p z^s\p\bar
z^s}+\frac{\p^2
h_{s\bar s}}{\p z^k\p\bar z^\ell}\right)\\
&+& \sum_{q,s}\left(\frac{\p h_{q\bar \ell}}{\p z^s}\frac{\p
h_{k\bar q}}{\p\bar z^s}+\frac{\p h_{k\bar q}}{\p z^s}\frac{\p
h_{q\bar \ell}}{\p\bar z^s}\right)-2\sum_{q,s}\left(\frac{\p
h_{q\bar \ell}}{\p z^k}\frac{\p h_{s\bar q}}{\p \bar z^s}+\frac{\p
h_{q\bar s}}{\p z^s}\frac{\p h_{k\bar q}}{\p \bar z^\ell}\right)\\
&=&\nonumber\frac{1}{2}\sum_s\left(\frac{\p^2 h_{s\bar \ell}}{\p
z^k\p\bar z^s}+\frac{\p^2 h_{k\bar s}}{\p z^s\p\bar z^\ell}\right)-
\sum_s\left(\frac{\p^2 h_{k\bar \ell}}{\p z^s\p\bar z^s}+\frac{\p^2
h_{s\bar s}}{\p z^k\p\bar z^\ell}\right)\\
&+& \frac{(T\circ \bar T)_{k\bar \ell}+(T\boxdot \bar T)_{k\bar
\ell}}{4}+\frac{1}{2}{(T(\p^*\omega)^\#)}_{i\bar
j}+\frac{1}{2}\bar{(T(\p^*\omega)^\#)}_{i\bar j}.\label{mathscrR}
\end{eqnarray}  \elemma

\vspace{0.5cm}

\subsection{Geometry of the first Levi-Civita Ricci curvature}

\bdefinition Let $M$ be a compact complex manifold. A Hermitian
metric $\omega$ on  $M$ is called \emph{balanced}, if $d^*\omega=0$.
$\omega$ is called \emph{conformally balanced}, if there exists a
smooth function $\phi:M\>\R$ and a balanced metric $\omega_B$ such
that $\omega=e^\phi\omega_B$. \edefinition

\btheorem \label{pluriclosed} On a compact Hermitian manifold
$(M,\omega)$, the first Levi-Civita Ricci form $\fR^{(1)}$
represents the first Aeppli-Chern class $ c_1^{AC}(M)$ in
$H^{1,1}_{A}(M)$. More precisely, \beq
\mathfrak{R}^{(1)}=\Theta^{(1)}-\frac{1}{2}(\p\p^*\omega+\bp\bp^*\omega).\label{key11}\eeq
   Moreover,

\bd
 \item  $\mathfrak{R}^{(1)}$ is
$d$-closed if and only if $\p\bp\bp^*\omega=0$;
\item if
$\bp\p^*\omega=0$, then $\mathfrak{R}^{(1)}$ represents the real
first Chern class $c_1(M)\in H^{2}_{dR}(M)$,  i.e.
$c_1(M)=c_1^{AC}(M)$ in $H^{2}_{dR}(M)$.;\item  if $\omega$ is
conformally balanced, then
 $\mathfrak{R}^{(1)}$ represents the first Chern class $c_1(M)\in
H^{1,1}_{\bp}(M)$ and also the first Bott-Chern class
$c^{BC}_1(M)\in H^{1,1}_{BC}(M)$;
\item $\mathfrak R^{(1)}=\Theta^{(1)}$ if and only if $d^*\omega=0$,
i.e. $(M,\omega)$ is a balanced manifold.\ed

 \etheorem

\bproof  By formulas (\ref{theta1}), (\ref{dd*o}), (\ref{pp*o}) and
(\ref{mathfrakR1}), we have \be &&\Theta^{(1)}_{i\bar
j}-\frac{1}{2}(\p\p^*\omega+\bp\bp^*\omega)_{i\bar
j}\\&=&-\sum_q\frac{\p^2 h_{q\bar q}}{\p z^i\p\bar
z^j}+\frac{1}{4}(T\boxdot\bar T)_{i\bar
j}-\frac{1}{2}\left(\sum_q\left(\frac{\p^2 h_{i\bar q}}{\p z^q\p\bar
z^j}+\frac{\p^2 h_{q\bar j}}{\p z^i\p\bar z^q}-2\frac{\p^2 h_{q\bar
q}}{\p z^i\p\bar z^j}\right) +(T\boxdot\bar T)_{i\bar
j}\right)\\
&=&-\frac{1}{2}\sum_q\left(\frac{\p^2 h_{i\bar q}}{\p z^q\p\bar
z^j}+\frac{\p^2 h_{q\bar j}}{\p z^i\p\bar
z^q}\right)-\frac{1}{4}(T\boxdot \bar T)_{i\bar j}
\\&=&\mathfrak{R}^{(1)}_{i\bar j}. \ee

\noindent By (\ref{key11}), we see $\mathfrak R^{(1)}$ represents
the first Aeppli-Chern class $ c_1^{AC}(M)$, i.e. $[\mathfrak
R^{(1)}]=[\Theta^{(1)}]$ as classes in $H^{1,1}_{A}(M)$.

Next we prove the  properties of $\mathfrak R^{(1)}$.

(1). By (\ref{key11}) again, $d\mathfrak
R^{(1)}=-\frac{1}{2}\left(\bp\p\p^*\omega+\p\bp\bp^*\omega\right)$.
By degree reasons, $d\mathfrak R^{(1)}=0$ if and only if
$\p\bp\bp^*\omega=0$.

(2).  If $\bp\p^*\omega=0$, we have
$$\mathfrak{R}^{(1)}=\Theta^{(1)}-\frac{1}{2}(\p\p^*\omega+\bp\bp^*\omega)=\Theta^{(1)}-\frac{1}{2}dd^*\omega.$$
Hence $\left[\mathfrak{R}^{(1)}\right]=\left[\Theta^{(1)}\right]\in
H^{2}_{dR}(M)$.

(3). If $\omega$ is conformally balanced, there exists a smooth
function $f$ and a balanced metric $ \omega_f$ such that
$\omega_f=e^f \omega$. We denote by an extra index $f$ the
corresponding quantities with respect to the new metric $\omega_f$.
 The Christoffel symbols of $\omega_f$ are $$(\Gamma_f)_{\bar i
j}^k=\frac{1}{2}e^{-f}g^{k\bar\ell}\left(\frac{\p (e^f g_{j\bar
\ell})}{\p\bar z^i}-\frac{\p (e^f g_{j\bar i})}{\p\bar
z^\ell}\right)=\Gamma_{\bar i
j}^k+\frac{1}{2}\left(\delta_{jk}f_{\bar i}-g^{k\bar \ell} g_{j\bar
i} f_{\bar \ell}\right).$$ In particular, $$ (\Gamma_f)_{\bar i
k}^k=\Gamma_{\bar i k}^k+\frac{n-1}{2}f_{\bar i}.$$ By (\ref{key1}),
we obtain \beq \bp^*_f\omega_f=\bp^*\omega+\sq(n-1)\p f.\eeq
Therefore, \beq  \bp\bp_f^*\omega_f=\bp\bp^*\omega-(n-1)\sq \p\bp f
\qtq{and} \p\bp_f^*\omega_f=\p\bp^*\omega.\label{conformal2}\eeq
Since $\omega_f$ is balanced, i.e. $\bp_f^*\omega_f=0$, we obtain
$$\p\p^*\omega+\bp\bp^*\omega=2(n-1)\sq\p\bp f.$$
Hence, $\mathfrak R^{(1)}=\Theta^{(1)}-(n-1)\sq\p\bp f$, i.e.
$\left[\mathfrak R^{(1)}\right]=\left[\Theta^{(1)}\right]\in
H^{1,1}_{BC}(M)$. Hence, $\fR^{(1)}$ represents the first Chern
class $c_1(M)\in H^{1,1}_{\bp}(M)$ and also the first Bott-Chern
class $c_1^{BC}(M)\in H^{1,1}_{BC}(M)$.

(4). $\mathfrak{R}^{(1)}=\Theta^{(1)}$ if and only if
$\p\p^*\omega+\bp\bp^*\omega=0$. By pairing with $\omega$, we see
the latter is equivalent to $d^*\omega=0$. \eproof

\bexample\label{ex} In this example, we shall construct a Hermitian
metric with strictly positive $\mathfrak{R}^{(1)}$, but
$\Theta^{(1)}$ is not strictly positive. Let $M$ be a Fano manifold
with complex dimension $n\geq 2$. By Yau's theorem, there exists a
K\"ahler metric $\omega$ on $M$ such that
$\mathfrak{R}^{(1)}_\omega=\Theta^{(1)}_\omega>0$. For any smooth
function $\phi$, we define $\omega_t=e^{t\phi}\omega$. By
(\ref{1chern}), one can see
\beq\Theta^{(1)}_{\omega_t}=\Theta^{(1)}_\omega-nt\sq\p\bp \phi.
\eeq Hence, by (\ref{key11}) and (\ref{conformal2}),\beq
\mathfrak{R}^{(1)}_{\omega_t}=\mathfrak R^{(1)}_\omega-\sq t \p\bp
\phi .\eeq
 Let $t_0=\sup\{t>0\ |\ \Theta^{(1)}_{\omega_t}=\Theta^{(1)}_\omega-nt\sq\p\bp\phi\geq 0\}$, and
$t_1:=\frac{3}{2} t_0$, then \beq
\mathfrak{R}^{(1)}_{\omega_{t_1}}=\fR^{(1)}_\omega-t_1\sq\p\bp\phi>0\eeq
but $\Theta^{(1)}_{\omega_{t_1}}$ is not positive definite.

\eexample

\bdefinition[\cite{DGMS75}] A compact \emph{complex} manifold $M$ is
said to satisfy the $\p\bp$-lemma if the following statement holds:
 if $\eta$ is $d$-exact, $\p$-closed and $\bp$-closed, it must be
 $\p\bp$-exact.
In particular, on such manifolds, for any pure-type form
$\psi\in\Om^{p,q}(M)$, if $\psi$ is $\bp$-closed and $\p$-exact,
then it is $\p\bp$-exact.\edefinition
 It is well-known  that all compact
K\"ahler manifolds satisfy the $\p\bp$-lemma. Moreover, if $\mu :
\hat M\>M$ is a modification between compact complex manifolds and
if the $\p\bp$-lemma holds for $\hat M$, then the $\p\bp$-lemma also
holds for $M$. In particular, Moishezon manifolds  and also
manifolds in Fujiki class $\mathscr C$  satisfy the $\p\bp$-lemma.
For more details, we refer to \cite{DGMS75,AT13,Po2,Po3} and also
the references therein.

 In the following, we show  on complex manifolds with $\p\bp$-lemma,
 the converse of $(2)$ and $(3)$ in Theorem \ref{pluriclosed} are
 also true.

\bproposition  Let $M$ be a compact complex manifold on which the
$\p\bp$-lemma holds. Let $\omega$ be a Hermitian metric on $M$. \bd

\item $\mathfrak{R}^{(1)}$ represents the real first Chern class $c_1(M)\in
H^{2}_{dR}(M)$  if and only if $\bp\p^*\omega=0$;
\item
 $\mathfrak{R}^{(1)}$ represents the first Chern class $c_1(M)\in
H^{1,1}_{\bp}(M)$ if and only if $\omega$ is conformally balanced.
 \ed
\eproposition \bproof(1) If $\left[\mathfrak
R^{(1)}\right]=\left[\Theta^{(1)}\right]\in H^{2}_{dR}(M)$, by
(\ref{key11}),
$\mathfrak{R}^{(1)}=\Theta^{(1)}-\frac{1}{2}dd^*\omega+\frac{1}{2}\left(\bp\p^*\omega+\p\bp^*\omega\right)$,
there exists a  $1$-form $\gamma$ such that
$d\gamma=\bp\p^*\omega+\p\bp^*\omega$. The compatibility condition
$d^2\gamma=0$ implies $\p\bp\p^*\omega+\bp\p\bp^*\omega=0$. Hence
$\bp\p\bp^*\omega=0$.
 If we set $
\gamma=\gamma^{1,0}+\gamma^{0,1}$, then
$\p\gamma^{1,0}=\p\bp^*\omega$. Therefore, $\p\gamma^{1,0}$ is both
$d$-closed and  $\p$-exact.  By $\p\bp$-lemma, there exists some
$\eta$ such that $\p\gamma^{1,0}=\p\bp\eta$. By degree reasons,
$\p\gamma^{1,0}=0$, that is $\p\bp^*\omega=0$.

\noindent (2). If $\left[\mathfrak
R^{(1)}\right]=\left[\Theta^{(1)}\right]\in H^{1,1}_{\bp}(M)$, there
exits $(1,0)$-form $\tau$ such that $\bp\tau=\p\p^*\omega$. So
$\bp\tau$ is both $\p$-closed and $\bp$-exact, hence by
$\p\bp$-lemma, there exists a smooth function $\phi$ such that
$\p\p^*\omega=\bp\tau=\sq\p\bp\phi$. By \cite{Ga3}, there exists a
smooth function $f$ such that $\omega_f:=e^{\frac{f}{n-1}}\omega$ is
a Gauduchon metric, i.e. $\p\bp\omega^{n-1}_f=0$. We use the index
$f$ to denote the operations with respect to the new metric
$\omega_f$. For example,
$$\|\p^*_f\omega_f\|^2_f=(\p\p^*_f\omega_f,\omega_f)_f.$$
We can see $\p\p^*_f\omega_f=\p\p^*\omega-\sq\p\bp
f=\sq\p\bp(\phi-f)$. Moreover,
$$(\p\p_f^*\omega_f, \omega_f)_f=(\sq\p\bp(\phi-f), \omega_f)_f=\int\sq\p\bp(\phi-f)\wedge \frac{\omega_f^{n-1}}{(n-1)!}=0,$$
since $\omega_f$ is a  Gauduchon metric. Therefore,
$\p_f^*\omega_f=0$, i.e. $\omega_f$ is balanced and so $\omega$ is
conformally balanced. \eproof

\bremark By (\ref{conformal2}), a conformally balanced metric
$\omega$ satisfies $\bp\p^*\omega=0$. On the other hand, if
$H_{\bp}^{0,1}(M)=0$, then $\bp\p^*\omega=0$ if and only if $\omega$
is conformally balanced. On the Hopf manifold $\S^{2n-1}\times \S^1$
with $n\geq 2$, the canonical metric $\omega_0$ satisfies
$\bp\p^*\omega_0=0$, but it is not conformally balanced.
 \eremark

\bcorollary\label{vanishingofchernclasses} Let $M$ be a complex
manifold. Then
$$c_1^{BC}(M)=0 \Longrightarrow c_1(M)=0 \Longrightarrow c_1^{AC}(M)=0.$$ Moreover, on a complex manifold
satisfying the $\p\bp$-lemma,
$$c_1^{BC}(M)=0\Longleftrightarrow c_1(M)=0\Longleftrightarrow
c_1^{AC}(M)=0.$$ \ecorollary

\bproof The first statement is obvious. For the second statement, we
only need to show that, on a complex manifold $M$ with $\p\bp$-lemma
if $c_1^{AC}(M)=0$, then $c_1^{BC}(M)=0$. Indeed, if
$c_1^{AC}(M)=0$, for a Hermitian metric $\omega$ on $M$,
$$\Theta^{(1)}=\p A+\bar{\p B}$$
where $A$ and $B$ are $(0,1)$ forms. It is obvious that $\p A$ is
$\bp$-closed and $\p$-exact, and so by $\p\bp$-lemma, there exists a
smooth function $f_1$ such that $\p A= \p\bp f_1$. Similarly, there
exists smooth function $f_2$ such that $\p B=\p\bp f_2$ and so
$\Theta^{(1)}=\p\bp (f_1-\bar{f_2})$. That is $c_1^{BC}(M)=0$.
\eproof

\bremark It is well-known that $\S^{2n-1}\times \S^1$ ($n\geq 2$)
has $c_1(M)=c_1^{AC}(M)=0$, but $c_1^{BC}(M)\neq 0$ (e.g.
\cite[Example ~3.3]{T}). It is also interesting to find a complex
manifold $N$ with $c^{AC}_1(N)=0$ but $c_1(N)\neq 0.$\eremark

\subsection{Hermitian manifolds with nonnegative  $\fR^{(1)}$} In
this subsection, we study the geometry of Hermitian manifolds with
nonnegative first  Levi-Civita Ricci curvature $\fR^{(1)}$.%
%

 \btheorem\label{main111}  Let $(M,h)$ be a compact
Hermitian manifold. If the first Levi-Civita Ricci curvature
$\fR^{(1)}$ is quasi-positive, then the top intersection number
$c^n_1(M)> 0$.  In particular, $H^2_{dR}(M)$, $H^{1,1}_{\bp}(M),
H^{1,1}_{BC}(M)$ and $H^{1,1}_A(M)$ are all non-zero.\etheorem

 \bproof At first, let's
recall the general theory for vector bundles. Let $\nabla^E$ be a
connection on the holomorphic vector bundle $E$. Let $r$ be the rank
of $E$, then there is a naturally induced connection
$\nabla^{\det(E)}$ on the determine line bundle
$\det(E)=\Lambda^rE$, \beq \nabla^{\det(E)}(s_1\wedge\cdots\wedge
s_r)=\sum_{i=1}^rs_1\wedge\cdots\wedge \nabla^E s_i\wedge
\cdots\wedge s_r.\eeq The curvature tensor of $(E,\nabla^E)$ is
denoted by $R^E\in\Gamma\left(M,\Lambda^2T^*M\ts End(E)\right)$ and
the curvature tensor of $(\det E, \nabla^{\det(E)})$ is denoted by
$R^{\det (E)}\in\Gamma\left(M,\Lambda^2T^*M\right).$ We have the
relation that \beq tr R^E=R^{\det
E}\in\Gamma(M,\Lambda^{2}T^*M).\eeq Note that the trace operator is
well-defined without using the metric on $E$. Moreover,
 $tr R^E=R^{\det E}$ is a $d$-closed
$2$-form.  By Bianchi identity, we know, for any vector bundle
$(F,\nabla^F)$
$$\nabla^{F\ts F^*} R^F=0.$$
In particular, if $F$ is a line bundle, $F\ts F^*=\underline{\C}$
and $\nabla^{F\ts F^*}=d$. Hence $d\left(R^{\det E}\right)=0$. On
the other hand, by Chern-Weil theory (e.g. \cite[Theorem~1.9]{Zh}),
$ R^{\det E}$ represents the real first Chern class $c_1(E) \in
H^2(M,\Z)$. In fact, let $\nabla^{Ch}$ be the Chern connection on
the Hermitian holomorphic line bundle $(\det E, h)$, and
$\Theta^{\det E}$ be the Chern curvature, then by Chern-Weil theory,
$$R^{\det E}-\Theta^{\det E}=d\beta$$ for some $1$-form $\beta$. It
is well known that the Chern curvature $\Theta^{\det E}$ of the
Hermitian line bundle $(\det E, h)$ represents the first Chern class
$c_1(E)\in H^{1,1}_{\bp}(M)$.

Now we go back to the setting on the Hermitian manifold
$(M,\omega)$. Let $E=T^{1,0}M$ with Hermitian metric $h$ induced by
$\omega$. With respect to the Levi-Civita connection $\hat \nabla$
on $E$, we have a decomposition $$ R^{\det
E}=\eta^{2,0}+\eta^{0,2}+\eta^{1,1}.\label{11}$$ It is obvious that
\beq \eta^{1,1}= \sq R_{i\bar j k}^kdz^i\wedge d\bar
z^j=\mathfrak{R}^{(1)},\ \ \
 \qtq{and}
\eta^{0,2}=\bar{\eta}^{2,0}.\eeq It is also easy to see that $
\eta^{2,0}=-\p\bp^*\omega$.
 Hence \beq\int
\left(R^{\det E}\right)^{n}=
\sum_{\ell=0}^{\left[\frac{n}{2}\right]}{n\choose
{2\ell}}{{2\ell}\choose \ell} \int(\eta^{2,0}\wedge \bar
\eta^{2,0})^{\ell}\wedge \left(\eta^{1,1}\right)^{n-2\ell}.\eeq It
is obvious that, if $\eta^{1,1}$ is quasi-positive, \beq
\int(\eta^{2,0}\wedge \bar \eta^{2,0})^{\ell}\wedge
\left(\eta^{1,1}\right)^{n-2\ell}\geq 0.\eeq for $1\leq \ell\leq
\left[\frac{n}{2}\right]$ and $\int \left(\eta^{1,1}\right)^{n}>0 $.
That is \beq \int \left(R^{\det E}\right)^{n}>0.\eeq We obtain \beq
\int c^n_1(M)= \int \left(R^{\det E}\right)^{n}>0.\eeq  On the other
hand, if $H^{1,1}_A(M)=0$, we obtain $\Theta^{(1)}=\p B+\bp C$ for
$1$-forms $B$ and $C$,  and so \beq \int_M
\left(\Theta^{(1)}\right)^n=\int_M (\p B+\bp C)\wedge
\left(\Theta^{(1)}\right)^{n-1}=0, \eeq which is a contradiction to
$c^n_1(M)>0$. The non-vanishing of other cohomology groups follows
immediately from Corollary \ref{vanishingofchernclasses}.
 \eproof
\bremark \bd\item Note that, in general, \beq \int_M
\left(\Theta^{(1)}\right)^n\neq \int_M\left(\mathfrak
R^{(1)}\right)^n.\eeq


\item When $M$ is in the Fujiki class $\mathscr C$ and $\fR^{(1)}$
is strictly positive, then $M$  is a K\"ahler manifold
(\cite[Theorem~0.2]{C}).

\ed

  \eremark


%
%
%
%

\subsection{Hypothetical complex structures on $\S^6$}
Let $(M,h)$ be a Hermitian manifold with constant Riemannian
sectional curvature $K$, i.e., for any $X,Y,Z,W\in T_\R M$, \beq
R(X,Y,Z,W)=K\cdot
\left(g(X,W)g(Y,Z)-g(X,Z)g(Y,W)\right).\label{constantse2}\eeq
Therefore, by the complexification process, $$ R_{i\bar j k\bar
\ell}=K\cdot h_{i\bar \ell}h_{k\bar j},\ \ R_{i\bar j}=K\cdot
h_{i\bar j}, \qtq{and} R_{ijk\bar \ell}=R_{\bar i\bar j k\bar
\ell}=0.$$ In particular, $Ric_H=K\cdot \omega.$ If $K>0$, we see
from (\ref{key0})  that \beq \mathfrak
R^{(1)}\geq Ric_H=K\cdot\omega_h>0.\eeq By Theorem \ref{main111}, $c_1^3(M)>0$. 
   Now we get Lebrun's
result that \bcorollary[\cite{Le87}] On $\S^6$, there is no
orthogonal complex structure compatible with metrics in  some small
neighborhood of the round metric. \ecorollary

\section{Curvature relations on Hermitian
manifolds}\label{riccirelations}
\subsection{Ricci curvature relations} Let's recall  different types of Ricci curvatures on a Hermitian manifold $(M,\omega)$:
\bd
\item the  Levi-Civita Ricci curvatures $\mathfrak{R}^{(1)}$, $\mathfrak{R}^{(2)}$, $\mathfrak{R}^{(3)}$, $\mathfrak{R}^{(4)}$;
\item the Chern Ricci curvatures $\Theta^{(1)}$, $\Theta^{(2)}$, $\Theta^{(3)}$,
$\Theta^{(4)}$;
\item the Hermitian-Ricci curvature $Ric_H=\sq R_{i\bar j} dz^i\wedge d\bar z^j$ where $R_{i\bar j}=h^{k\bar \ell} R_{i\bar j k\bar
\ell}$ (which is equal to $R^{(1)}$ and $R^{(2)}$), the third and
fourth Hermitian-Ricci curvatures $R^{(3)}$ and $R^{(4)}$;
\item the $(1,1)$-component of the complexified
Riemannian Ricci curvature,  $\mathscr Ric$.
 \ed In this subsection, we shall explore explicit relations between
them by using  $\omega$ and its torsion $T$.  We shall prove Theorem
\ref{riccire}. We also state it as in the following to the reader's
convenience.

\btheorem\label{main-2} Let $(M,\omega)$ be a compact Hermitian
manifold. \bd

\item The  Levi-Civita Ricci curvatures are  \beq
\mathfrak{R}^{(1)}=\Theta^{(1)}-\frac{1}{2}\left(\p\p^*\omega+\bp\bp^*\omega\right);\label{53}\eeq
\beq \mathfrak{R}^{(2)}=\Theta^{(1)}-\frac{1}{2}\left(
\p\p^*\omega+\bp\bp^*\omega\right)-\frac{\sq}{4}T\circ\bar
T+\frac{\sq}{4} T\boxdot \bar T;\label{54}\eeq \beq
\mathfrak{R}^{(3)}=\Theta^{(1)}-\frac{1}{2}\left(\sq\Lambda\left(\p\bp\omega\right)+(\p\p^*\omega+\bp\bp^*\omega)\right)+\frac{\sq}{4}
T\boxdot \bar T+\frac{T([\p^*\omega]^{\#})}{4};\label{553}\eeq \beq
\mathfrak{R}^{(4)}=\Theta^{(1)}-\frac{1}{2}\left(\sq\Lambda\left(\p\bp\omega\right)+(\p\p^*\omega+\bp\bp^*\omega)\right)+\frac{\sq}{4}
T\boxdot \bar T+\frac{\bar{T([\p^*\omega]^{\#})}}{4},\label{554}\eeq
where $(\p^*\omega)^{\#}$ is the dual vector of  the $(0,1)$-form
$\p^*\omega$.

\item  The  Chern-Ricci curvatures are \beq
\Theta^{(2)}=\Theta^{(1)}-\sq\Lambda\left(\p\bp\omega\right)-(\p\p^*\omega+\bp\bp^*\omega)+{\sq}T\boxdot\bar
T;\label{563}\eeq

\beq \Theta^{(3)}=\Theta^{(1)}-\p\p^*\omega;\label{564}\eeq

\beq \Theta^{(4)}=\Theta^{(1)}-\bp\bp^*\omega.\label{565}\eeq

\item The Hermitian-Ricci curvatures are \beq  Ric_H=R^{(1)}=R^{(2)}=\Theta^{(1)}-\frac{1}{2}\left(
\p\p^*\omega+\bp\bp^*\omega\right)-\frac{\sq}{4}T\circ\bar
T;\label{key2}\eeq \begin{eqnarray}
R^{(3)}=R^{(4)}\nonumber&=&\Theta^{(1)}-\frac{1}{2}\left(\sq\Lambda\left(\p\bp\omega\right)+(\p\p^*\omega+\bp\bp^*\omega)\right)\\&&+\frac{\sq}{4}
T\boxdot \bar
T+\frac{\bar{T([\p^*\omega]^{\#})}+T([\p^*\omega]^{\#})}{4}.\label{653}\end{eqnarray}

\item The $(1,1)$-component of the Riemannian Ricci
curvature is \begin{eqnarray} \mathscr
{R}ic\nonumber&=&\Theta^{(1)}-\sq(\Lambda
\p\bp\omega)-\frac{1}{2}(\p\p^*\omega+\bp\bp^*\omega)+\frac{\sq}{4}\left(2T\boxdot
\bar T+T\circ \bar T\right)\\&&+\frac{1}{2}\left(
T([\p^*\omega]^{\#})+\bar{
T([\p^*\omega]^{\#})}\right).\label{753}\end{eqnarray}\ed
 \etheorem

\bproof (1). Equation (\ref{53}) is proved in (\ref{key11}). By
(\ref{mathfrakR1}), (\ref{mathfrakR2}) and (\ref{53}), we see \be
\mathfrak R^{(2)}&=&\mathfrak R^{(1)}-\frac{\sq}{4}T\circ\bar
T+\frac{\sq}{4} T\boxdot \bar T\\&=&\Theta^{(1)}-\frac{1}{2}\left(
\p\p^*\omega+\bp\bp^*\omega\right)-\frac{\sq}{4}T\circ\bar
T+\frac{\sq}{4} T\boxdot \bar T\ee which proves (\ref{54}). A
straightforward computation shows \beq \sq
\Lambda(\p\bp\omega)=\sq\sum_q \left[\left(\frac{\p h_{i\bar j}}{\p
z^q \p\bar z^q}+\frac{\p h_{q\bar q}}{\p z^i\p\bar
z^j}\right)-\left(\frac{\p h_{i\bar q}}{\p z^q\p\bar z^j}+\frac{\p^2
h_{q\bar j}}{\p z^i\p\bar z^q}\right)\right]dz^i\wedge d\bar
z^j.\label{76}\eeq  If we write
$$
\Theta^{(1)}-\frac{1}{2}\left(\sq\Lambda\left(\p\bp\omega\right)+(\p\p^*\omega+\bp\bp^*\omega)\right)+\frac{\sq}{4}
T\boxdot \bar T:=\sq B_{i\bar j}dz^i\wedge d\bar z^j,$$ by using
formulas (\ref{theta1}), (\ref{76}), (\ref{dd*o}) and (\ref{pp*o}),
we obtain \be B_{i\bar j}&=&-\sum_q\frac{\p^2 h_{q\bar q}}{\p
z^i\p\bar z^j}-\frac{1}{2}\sum_q \left[\left(\frac{\p h_{i\bar
j}}{\p z^q \p\bar z^q}+\frac{\p h_{q\bar q}}{\p z^i\p\bar
z^j}\right)-\left(\frac{\p h_{i\bar q}}{\p z^q\p\bar z^j}+\frac{\p^2
h_{q\bar j}}{\p z^i\p\bar
z^q}\right)\right]\\&&-\frac{1}{2}\sum_q\left(\frac{\p^2 h_{i\bar
q}}{\p z^q\p\bar z^j}+\frac{\p^2 h_{q\bar j}}{\p z^i\p\bar
z^q}-2\frac{\p^2 h_{q\bar
q}}{\p z^i\p\bar z^j}\right) \\
&=&-\frac{1}{2}\sum_q\left(\frac{\p^2 h_{q\bar q}}{\p z^i\p\bar
z^j}+\frac{\p^2 h_{i\bar j}}{\p z^q\p\bar z^q}\right).\ee By formula
(\ref{mathfrakR3}), we obtain (\ref{553}). The proof of formula
(\ref{554}) follows from (\ref{mathfrakR4}), (\ref{mathfrakR3}) and
(\ref{553}).

(2). For Chern-Ricci curvature relations, we use similar
computations. For example, if we write
$$
\Theta^{(1)}-\left(\sq\Lambda\left(\p\bp\omega\right)+(\p\p^*\omega+\bp\bp^*\omega)\right)+{\sq}
T\boxdot \bar T:=\sq C_{i\bar j}dz^i\wedge d\bar z^j,$$ by using
formulas (\ref{theta1}), (\ref{76}), (\ref{dd*o}) and (\ref{pp*o}),
we obtain \be C_{i\bar j}&=&-\sum_q\frac{\p^2 h_{q\bar q}}{\p
z^i\p\bar z^j}-\sum_q \left[\left(\frac{\p h_{i\bar j}}{\p z^q
\p\bar z^q}+\frac{\p h_{q\bar q}}{\p z^i\p\bar
z^j}\right)-\left(\frac{\p h_{i\bar q}}{\p z^q\p\bar z^j}+\frac{\p^2
h_{q\bar j}}{\p z^i\p\bar
z^q}\right)\right]\\&&-\sum_q\left(\frac{\p^2 h_{i\bar q}}{\p
z^q\p\bar z^j}+\frac{\p^2 h_{q\bar j}}{\p z^i\p\bar z^q}-2\frac{\p^2
h_{q\bar
q}}{\p z^i\p\bar z^j}\right)+\frac{1}{4}\left(T\boxdot\bar T\right)_{i\bar j} \\
&=&-\sum_q\frac{\p^2 h_{i\bar j}}{\p z^q\p\bar
z^q}+\frac{1}{4}\left(T\boxdot\bar T\right)_{i\bar j}.\ee By formula
(\ref{theta2}), we obtain (\ref{563}). Formula (\ref{564}) follows
from (\ref{theta1}), (\ref{theta3}) and (\ref{dd*o}). Formula
(\ref{565}) follows from (\ref{theta1}), (\ref{theta4}) and
(\ref{pp*o}).

(3). For Hermitian-Ricci curvatures, by (\ref{R1R2}) and
(\ref{mathfrakR1}), we see
$$R^{(1)}=R^{(2)}=\mathfrak R^{(1)}-\frac{\sq}{4}T\circ \bar T.$$
Therefore, (\ref{key2}) follows from (\ref{53}). Similarly, by
(\ref{R3R4}), (\ref{mathfrakR3}) and (\ref{553}), we obtain
(\ref{653}).

(4). If we write
$$
\Theta^{(1)}-\left(\sq\Lambda\left(\p\bp\omega\right)+\frac{1}{2}(\p\p^*\omega+\bp\bp^*\omega)\right):=\sq
F_{i\bar j}dz^i\wedge d\bar z^j,$$ by using formulas (\ref{theta1}),
(\ref{76}), (\ref{dd*o}) and (\ref{pp*o}), we obtain \be F_{i\bar
j}&=&-\sum_q\frac{\p^2 h_{q\bar q}}{\p z^i\p\bar z^j}-\sum_q
\left[\left(\frac{\p h_{i\bar j}}{\p z^q \p\bar z^q}+\frac{\p
h_{q\bar q}}{\p z^i\p\bar z^j}\right)-\left(\frac{\p h_{i\bar q}}{\p
z^q\p\bar z^j}+\frac{\p^2 h_{q\bar j}}{\p z^i\p\bar
z^q}\right)\right]\\&&-\frac{1}{2}\sum_q\left(\frac{\p^2 h_{i\bar
q}}{\p z^q\p\bar z^j}+\frac{\p^2 h_{q\bar j}}{\p z^i\p\bar
z^q}-2\frac{\p^2 h_{q\bar
q}}{\p z^i\p\bar z^j}\right)-\frac{1}{4}\left(T\boxdot\bar T\right)_{i\bar j} \\
&=&\frac{1}{2}\sum_q\left(\frac{\p^2 h_{q\bar j}}{\p z^i\p\bar
z^q}+\frac{\p^2 h_{i\bar q}}{\p z^q\p\bar z^j}\right)-
\sum_q\left(\frac{\p^2 h_{i\bar j}}{\p z^q\p\bar z^q}+\frac{\p^2
h_{q\bar q}}{\p z^i\p\bar z^j}\right)-\frac{1}{4}\left(T\boxdot\bar
T\right)_{i\bar j}.\ee By (\ref{mathscrR}), we obtain (\ref{753}).
 \eproof

\subsection{Scalar curvature relations}
On a Hermitian manifold $(M,\omega)$, we can define five different
types of scalar curvatures:\bd
\item $s$, the scalar curvature of the background Riemannian metric;\item $\displaystyle s_R=h^{i\bar l}h^{k\bar j}R_{i\bar j k\bar
\ell}$, the Riemannian type scalar curvature;
\item $s_H=h^{i\bar j}h^{k\bar \ell}R_{i\bar j k\bar \ell}$, the scalar curvature of
the Hermitian curvature;
\item $s_{LC}=h^{i\bar j k\bar \ell}\mathfrak{R}_{i\bar j k\bar \ell}$, the scalar curvature of the
Levi-Civita connection;\item $s_C=h^{i\bar j} h^{k\bar \ell}
\Theta_{i\bar j k\bar \ell}$, the scalar curvature of the Chern
connection. \ed

\noindent By Theorem \ref{main-2}, we get the  corresponding scalar
curvature relations: \bcorollary\label{A} Let $(M,\omega)$ be a
compact Hermitian manifold, then \beq
s=2s_C+\left(\la\p\p^*\omega+\bp\bp^*\omega,\omega\ra-2|\p^*\omega|^2\right)-\frac{1}{2}|T|^2,\label{key-1}\eeq
\beq
s_{LC}=s_C-\frac{1}{2}\la\p\p^*\omega+\bp\bp^*\omega,\omega\ra=s_C-\la\p\p^*\omega,\omega\ra,\label{010}\eeq
\beq
s_H=s_C-\frac{1}{2}\la\p\p^*\omega+\bp\bp^*\omega,\omega\ra-\frac{|T|^2}{4},\label{009}\eeq
and \beq
s_R=s_C-\frac{1}{2}|\p^*\omega|^2-\frac{1}{4}|T|^2.\label{008}\eeq

 \ecorollary

\bproof By (\ref{scalar2}), we know $s=2h^{i\bar j}\mathscr R_{i\bar
j}$. By (\ref{753}),
$$s=2s_C-\left\la2\sq(\Lambda
\p\bp\omega),\omega\right\ra-\left\la\p\p^*\omega+\bp\bp^*\omega,\omega\right\ra+\frac{3|T|^2}{2}-2|\p^*\omega|^2,$$
where we use the fact that $ tr_\omega
T([\p^*\omega]^{\#})=-|\p^*\omega|^2$(see (\ref{98})).  By
(\ref{563}), we have \beq
|T|^2=\la\sq\Lambda\p\bp\omega,\omega\ra+\la\p\p^*\omega+\bp\bp^*\omega,\omega\ra,\label{key-3}
\eeq we get (\ref{key-1}).  (\ref{010}) and (\ref{009}) follow from
 (\ref{53}) and (\ref{key2}) respectively. For (\ref{008}), it
 follows from (\ref{key-1}),
 (\ref{010}) and the fact that $s=4s_R-2s_H$.
\eproof

\bremark\label{remark3}  Note  that we can also define the scalar
curvatures as
$$\tr_\omega \mathfrak R^{(3)},\ \ \  \tr_\omega\Theta^{(3)} \qtq{and} \tr_\omega R^{(3)}$$
The corresponding scalar curvature relations follow immediately from
formulas (\ref{553}), (\ref{564}) and (\ref{653}). More precisely,
we obtain \beq \tr_\omega \mathfrak R^{(3)}=\tr_\omega \mathfrak
R^{(4)}=s_C-\frac{1}{4}|T|^2-\frac{1}{4}|\p^*\omega|^2;\eeq \beq
\tr_\omega\Theta^{(3)}=\tr_\omega\Theta^{(4)}=s_C-
\la\p\p^*\omega,\omega\ra=s_{LC};\eeq  and \beq \tr_\omega
R^{(3)}=\tr_\omega
R^{(4)}=s_C-\frac{1}{4}|T|^2-\frac{1}{2}|\p^*\omega|^2=s_R.\eeq
 \eremark

\bcorollary\label{B} Let $(M,\omega)$ be a compact Hermitian
manifold. Then the following are equivalent: \bd

\item $(M,\omega)$ is K\"ahler;

\item $\displaystyle \int s\cdot \omega^n=\int 2s_C\cdot \omega^n$;

\item $\displaystyle \int s_C\cdot \omega^n=\int s_R\cdot \omega^n$;

\item $\displaystyle \int s_C\cdot \omega^n=\int s_{H}\cdot \omega^n$;

\item $\displaystyle \int s_H\cdot \omega^n=\int s_{LC}\cdot
\omega^n$.\ed \ecorollary

\bcorollary\label{C} Let $(M,\omega)$ be a compact Hermitian
manifold. Then the following are equivalent: \bd

\item $(M,\omega)$ is balanced;

\item $\displaystyle \int s\cdot \omega^n=\int 2s_R\cdot \omega^n$;

\item $\displaystyle \int s\cdot \omega^n=\int 2s_H\cdot \omega^n$;

\item $\displaystyle \int s_C\cdot \omega^n=\int s_{LC}\cdot \omega^n$;

\item $\displaystyle \int s_R\cdot \omega^n=\int s_H\cdot \omega^n$.

\ed\ecorollary

\section{Special metrics on Hermitian manifolds}\label{special}

 Before discussing special metrics on Hermitian manifolds,
we need the following observation which is the integral version of
(\ref{key-3}). We assume $\dim_\C M=n\geq 3$.

\bproposition\label{key-6} On a compact Hermitian manifold
$(M,\omega)$, for any $1\leq k\leq n-1$, we have \beq  \int \sq
\p\omega\wedge \bp\omega \wedge
\frac{\omega^{n-3}}{(n-3)!}=\|\p^*\omega\|^2-\|\p\omega\|^2,\eeq and
\beq \int \sq \omega^{n-k-1}\wedge \p\bp\omega^{k}=
(n-3)!k(n-k-1)\left(\|\p\omega\|^2-\|\p^*\omega\|^2\right).\label{key-4}\eeq

\eproposition

\bproof At first, it is easy to see \beq
\p\bp\omega^k=k\omega^{k-1}\p\bp\omega+k(k-1)\omega^{k-2}\p\omega\wedge
\bp\omega.\label{88}\eeq On the other hand, the $(1,2)$-form
$\alpha:=\bp\omega-\frac{L\Lambda\bp\omega}{n-1}$ is primitive, i.e.
$\Lambda\alpha=0$. Hence, by \cite[Proposition~6.29]{Voisin02}, \beq
*\left(\alpha\right)=(-1)^{\frac{3(3+1)}{2}} (\sq)^{1-2}\frac{L^{n-3}}{(n-3)!} \alpha =-\sq \frac{\omega^{n-3}\wedge \alpha}{(n-3)!}.\eeq
Therefore \beq \left(\p\omega,
\p\omega-\frac{L\Lambda\p\omega}{n-1}\right)=(\p\omega,\bar
\alpha)=\int \p\omega\wedge
*(\alpha)=-\int \sq \p\omega\wedge \alpha \wedge \frac{\omega^{n-3}}{(n-3)!}.\eeq
In particular, we have \be \int \sq \p\omega\wedge \bp\omega \wedge
\frac{\omega^{n-3}}{(n-3)!}&=&-\left(\p\omega,
\p\omega-\frac{L\Lambda\p\omega}{n-1}\right)+\int \sq \p\omega\wedge
\frac{L\Lambda\bp\omega}{n-1} \wedge \frac{\omega^{n-3}}{(n-3)!}\\
&=&-\|\p\omega\|^2+\frac{\|\Lambda\p\omega\|^2}{n-1}+\int \sq
\p\omega\wedge \frac{\Lambda\bp\omega}{n-1} \wedge
\frac{\omega^{n-2}}{(n-3)!}. \ee Form (\ref{key1}),
$\p^*\omega=-\sq\Lambda\bp\omega$, we have \be &&\int \sq
\p\omega\wedge \frac{\Lambda\bp\omega}{n-1} \wedge
\frac{\omega^{n-2}}{(n-3)!}\\&=&\frac{1}{(n-1)^2}\int
\p^*\omega\wedge \frac{\p
\omega^{n-1}}{(n-3)!}=\frac{n-2}{n-1}\int\p\p^*\omega\wedge
\frac{\omega^{n-1}}{(n-1)!}\\
&=&\frac{n-2}{n-1}(\p\p^*\omega,\omega)=\frac{n-2}{n-1}\|\p^*\omega\|^2.
\ee Therefore \beq  \int \sq \p\omega\wedge \bp\omega \wedge
\frac{\omega^{n-3}}{(n-3)!}=\|\p^*\omega\|^2-\|\p\omega\|^2.\label{89}\eeq
From integration by parts, we also get \beq \int
\sq\p\bp\omega\wedge
\frac{\omega^{n-2}}{(n-2)!}=\|\p\omega\|^2-\|\p^*\omega\|^2.\label{999}\eeq
Hence (\ref{key-4})  follows from (\ref{88}), (\ref{89}) and (\ref{999}).%
%
%
%
%
%
 \eproof

\bcorollary We have the follow relation on a compact Hermitian
manifold $(M,\omega)$ \beq \int s_R\cdot \frac{\omega^n}{n!}=\int
s_{LC}\cdot \frac{\omega^n}{n!}-\frac{1}{(n-3)!2k(n-k-1)}\int \sq
\omega^{n-k-1}\wedge \p\bp\omega^k.\eeq \bproof It follows from
Corollary \ref{A} and formula (\ref{key-4}). \eproof \ecorollary

\noindent  Fu-Wang-Wu defined in \cite{FWW10} that a Hermitian
metric $\omega$ satisfying \beq \sq \omega^{n-k-1}\wedge
\p\bp\omega^{k}=0,\ \ \ 1\leq k\leq n-1\eeq is called a
$k$-Gauduchon metric. It is obvious that $(n-1)$-Gauduchon metric,
i.e. $\p\bp\omega^{n-1}=0$ is the original Gauduchon metric. It is
well-known that, Hopf manifolds $\S^{2n+1}\times S^1$ can not
support Hermitian metrics with $\p\bp\omega=0$ (SKT) or
$d^*\omega=0$ (balanced). They showed in \cite{FWW10} that on
$\S^5\times \S^1$, there exists a $1$-Gauduchon metric $\omega$,
i.e. $\omega\wedge \p\bp\omega=0$.

 As a straightforward application of Proposition \ref{key-6}, we obtain:
 \bcorollary If $(M,\omega)$ is $k$-Gauduchon ($1\leq k\leq n-2$) and
 also balanced, then $(M,\omega)$ is K\"ahler.
 \ecorollary
This is also true in the ``conformal" setting:

\bcorollary On a compact complex manifold, the following are
equivalent:

\bd \item $(M,\omega)$ is conformally K\"ahler;

\item $(M,\omega)$ is conformally $k$-Gauduchon for $1\leq k\leq
n-2$, and conformally balanced;

\ed

In particular, the following are also equivalent: \bd \item[(3)]
$(M,\omega)$ is K\"ahler;

\item[(4)] $(M,\omega)$ is  $k$-Gauduchon for $1\leq k\leq
n-2$, and conformally balanced;

\item[(5)]   $(M,\omega)$ is conformally balanced and $\Lambda^2(\p\bp\omega)=0$.
\ed \ecorollary

\bproof We first show $(2)$ implies $(1)$. Since $\omega$ is
conformally balanced, $\omega=e^F\omega_B$ for a balanced metric
$\omega_B$ and a smooth function $F\in C^\infty(M,\R)$.  By the
 conformally $k$-Gauduchon condition, we know there exists $\tilde F\in
 C^{\infty}(M,\R)$ and a $k$-Gauduchon metric $\omega_G$ such that $\omega=e^{\tilde
 F}\omega_G$. Let $f=F-\tilde F$, then $\omega_G=e^f\omega_B$.
 Since $\omega_G$ is $k$-Gauduchon,
$$\left(e^{f}\omega_B\right)^{n-k-1}\wedge
\p\bp\left(e^f\omega_B\right)^k=0,$$ and we obtain \beq
\omega^{n-k-1}_B\wedge \p\bp (e^f\omega_B)^k=0.\label{1010}\eeq

\noindent\textbf{Claim:} If a balanced metric $\omega_B$ satisfies
(\ref{1010}), then $f$ is a constant and $\omega_B$ is K\"ahler.\\

Since $\omega_B$ is balanced, i.e.,
$\p\omega_B^{n-1}=\bp\omega_B^{n-1}=0$, we see
$\omega_B^{n-k-1}\wedge \p\omega^k_B=0$. From (\ref{1010}), we get
\beq e^{kf}\omega_B^{n-k-1}\wedge
\p\bp\omega_B^k+\omega^{n-1}_B\wedge \p\bp
\left({e^{kf}}\right)=0.\label{101}\eeq Hence, \beq
\int_Me^{kf}\cdot\omega_B^{n-k-1}\wedge
\p\bp\omega_B^k=-\int_M\omega^{n-1}_B\wedge \p\bp
\left({e^{kf}}\right)=0. \eeq Using integration by parts and the
balanced condition  $\omega_B^{n-k-1}\wedge \p\omega_B^{k}=0$, we
obtain \be 0&=& -\sq\int_Me^{kf}\cdot\omega_B^{n-k-1}\wedge
\p\bp\omega_B^k\\&=&\sq k(n-k-1) \int_M e^{kf} \omega_B^{n-3}
\p\omega_B \wedge \bp\omega_B\\&=&\sq k(n-k-1)(n-3)!\int_M
e^{\frac{kf}{2}}\p\omega_B\wedge
\left(e^{\frac{kf}{2}}\bp\omega_B\cdot
\frac{\omega^{n-3}_B}{(n-3)!}\right)
\\
&=&-k(n-k-1)(n-3)!\|e^{\frac{kf}{2}}\p\omega_B\|^2_B, \ee and so
$\p\omega_B=0$, i.e. $\omega_B$ is  K\"ahler. Note that, in the last
step we use the fact that $e^{\frac{kf}{2}}\bp\omega_B$ is
primitive, i.e.
$$\Lambda\left(e^{\frac{kf}{2}}\bp\omega_B\right)=e^{\frac{kf}{2}}\Lambda\bp\omega_B=\sq e^{\frac{kf}{2}}\p^*\omega_B=0,$$
where the norm $\|\bullet\|_B$, $\p^*$ and the contraction $\Lambda$
are taken with respect to $\omega_B$. From (\ref{101}), we see that
$f$ must be constant if $\omega_B$ is K\"ahler. The proof of the
claim is complete. We know $\omega$ is conformally K\"ahler.

The equivalence of $(4)$ and $(3)$ follows from the proof of the
claim in the last paragraph under the condition $\omega=\omega_G$,
i.e. $\tilde F=0$ and $f=F$. Next we show $(5)$ implies $(3)$. In
fact, if $\omega$ is conformally balanced, i.e. $\omega=e^f\omega_B$
for some balanced metric $\omega_B$ and smooth function $f\in
C^\infty(M,\R)$, the condition $\Lambda^2\p\bp\omega=0$
 implies $$\Lambda_B^2\p\bp\left(e^f\omega_B\right)=0$$ where
 $\Lambda_B$ is the contraction operator with respect to $\omega_B$.
 By duality, we have
 \be 0=\int_M \p\bp(e^f \omega_B)\wedge \omega^{n-2}_B=\int_M
 e^f\omega_B\wedge \p\bp \omega_B^{n-2}.
 \ee
As similar as the proof in the last paragraph, we obtain both
$\omega$ and $\omega_B$ are K\"ahler. \eproof

\section{Levi-Civita Ricci-flat  and constant negative scalar curvature metrics on Hopf manifolds}\label{example}

In this section, we construct special Hermitian metrics on
non-K\"ahler manifolds related to Hopf manifolds. More precisely,
\bd \item We construct explicit Levi-Civita Ricci-flat metrics on
$\S^{2n-1}\times \S^1$;
\item We construct a smooth family of explicit Hermitian metrics $h_\lambda$ with
$\lambda\in(-1,+\infty)$ on $\S^{2n-1}\times \S^1$ such that their
Riemannian scalar curvature are constants and vary from a positive
constant to $-\infty$. In particular, we obtain Hermitian metrics
with negative constant Riemannian scalar curvature on Hermitian
manifolds with $c_1\geq 0$;

\item We construct pluriclosed metrics on the projective bundles over
$\S^{2n-1}\times \S^1$. \ed

\subsection{Levi-Civita Ricci-flat metrics on Hopf manifolds}
 Let's recall an example in
\cite[Section~6]{LY12}. Let $M=\S^{2n-1}\times \S^1$ be the standard
$n$-dimensional ($n\geq 2$) Hopf manifold. It is diffeomorphic to
$\C^n- \{0\}/G$ where $G$ is cyclic group generated by the
transformation $z\rightarrow \frac{1}{2}z$. It has an induced
complex structure from $\C^n-\{0\}$.   On $M$, there is a natural
induced metric $\omega_0$ given by
 \beq \omega_0=\sq h_{i\bar j}dz^i\wedge d\bar z^j=
\frac{4\delta_{i\bar j}}{|z|^2}dz^i\wedge  d\bar z^j.\label{hopf}
\eeq For the reader's convenience, we include some elementary
computations:
 $$\frac{\p h_{k\bar
\ell}}{\p z^i}=-\frac{4\delta_{k\ell}\bar z^i}{|z|^4},\ \ \ \
\frac{\p h_{k\bar \ell}}{\p\bar z^j}=-\frac{4
\delta_{k\ell}z^j}{|z|^4} $$ and $$ \frac{\p^2 h_{k\bar \ell}}{\p
z^i\p\bar z^j}=-4\delta_{k\ell}\frac{\delta_{i\bar j}|z|^2-2\bar z^i
z^j}{|z|^6}.$$
 Similarly, we have
$$ \Gamma_{ik}^\ell =\frac{1}{2}h^{\ell\bar q}\left(\frac{\p h_{i\bar
q}}{\p z^k}+\frac{\p h_{k\bar q}}{\p
z^i}\right)=-\frac{\delta_{i\ell}\bar z^k+\delta_{k\ell}\bar
z^i}{2|z|^2},$$$$ \Gamma_{\bar j k}^\ell=\frac{1}{2}h^{\ell\bar
q}\left(\frac{\p h_{k\bar q}}{\p \bar z^j}-\frac{\p h_{k\bar j}}{\p
\bar
z^q}\right)=\frac{\delta_{jk}z^\ell-\delta_{k\ell}z^j}{2|z|^2},$$
and
$$ \frac{\p\Gamma_{ik}^\ell}{\p \bar
z^j}=-\frac{\delta_{{k\ell}}\delta_{ij}+\delta_{i\ell}\delta_{jk}}{2|z|^2}+\frac{\delta_{i\ell}z^j\bar
z^k+\delta_{k\ell}z^j\bar z^i}{2|z|^4}, $$ $$ \frac{\p \Gamma_{\bar
j k}^\ell}{\p
z^i}=\frac{\delta_{jk}\delta_{i\ell}-\delta_{k\ell}\delta_{ij}}{2|z|^2}-\frac{(\delta_{jk}z^\ell-\delta_{k\ell}z^j)\bar
z^i}{2|z|^4}.$$
 The  curvature components of $\mathfrak{R}$ are \beq\mathfrak{R}_{i\bar j
k\bar\ell}=-h_{p\bar\ell}\left(\frac{\p \Gamma^{p}_{ik}}{\p \bar
z^j}-\frac{\p \Gamma^{p}_{\bar jk}}{\p z^i}+\Gamma_{
ik}^{s}\Gamma^{p}_{\bar js}-\Gamma_{ \bar jk}^{s}\Gamma^{p}_{s
i}\right)=\frac{3\delta_{i\ell}\delta_{jk}}{|z|^4}-\frac{2\delta_{i\ell}
z^j\bar z^k+2\delta_{jk}z^\ell\bar z^i-\delta_{ij}\bar
z^kz^\ell}{|z|^6}.\eeq
 The complexified  curvature components are  \beq R_{i\bar j
k\bar\ell}=\mathfrak{R}_{i\bar j
k\bar\ell}+h_{p\bar\ell}\Gamma_{\bar s i}^p\bar{\Gamma_{\bar k
j}^s}= \frac{2\delta_{i\ell}\delta_{jk}}{|z|^4}-\frac{\delta_{i\ell}
z^j\bar z^k+\delta_{jk}z^\ell\bar z^i}{|z|^6}.\eeq The Chern
curvature components are \beq \Theta_{i\bar j k\bar
\ell}=-\frac{\p^2 h_{k\bar \ell}}{\p z^i\p\bar z^j}+h^{p\bar
q}\frac{\p h_{{k\bar q}}}{\p z^i}\frac{\p h_{p\bar \ell}}{\p\bar
z^j}=\frac{4\delta_{kl}(\delta_{ij}|z|^2-z^j\bar z^i)}{|z|^6}.
\label{Griffiths}\eeq

\noindent  As  consequences (see also \cite{LY12}), \beq
\Theta^{(1)}=-\sq\p\bp\log\det(h)=n\cdot\sq \p\bp\log |z|^2,\ \
\Theta^{(2)}=\frac{n-1}{4}\omega_h;\label{dd}\eeq \beq \mathfrak
R^{(1)}=\sq\p\bp\log |z|^2,\ \ \ \
\mathfrak{R}^{(2)}=\frac{4-n}{4}\cdot \sq \p\bp\log
|z|^2+\frac{n-1}{16}\omega_h;\eeq \beq Ric_H=\frac{\sq}{2}\p\bp\log
|z|^2;\label{1000001}\eeq \beq \mathscr Ric=\frac{n-1}{2}\cdot \sq
\p\bp\log |z|^2+\frac{n-1}{8}\omega_h\eeq and \beq
s_R=\frac{n^2-n}{8},\ \ s_H=\frac{n-1}{8}, \ \
s=\frac{(2n-1)(n-1)}{4}, \ \ s_{LC}=\frac{n-1}{4}, \ \ \
s_{C}=\frac{n(n-1)}{4}.\eeq

\bremark\label{remark2} Let $M=\S^{2n-1}\times \S^1$ with $n\geq 2$.

 \bd \item Since $\p\log|z|^2$ is a globally defined $(1,0)$ form on $M$, one can see from $(\ref{dd})$ that $c_1(M)=0$, $c^{AC}_1(M)=0$, but $c_1^{BC}(M)\neq 0$. In particular, the canonical line
 bundle is topologically trivial but not holomorphically trivial;

\item Since the Chern scalar curvature $s_C=\frac{n(n-1)}{4}$, one can see that, $H^0(M,mK_M)=0$ for any integer
$m\geq 1$. In particular, the Kodaira dimension of  $M$ is
$-\infty$, and so $K_M$ is not a torsion line bundle;

\item The canonical metric $\omega_h$ satisfies $\bp\p^*\omega_h=0$, but it is well-known that $M$ can not support any balanced
metric;

\item $\p\bp$-lemma does not hold on $M$;

 \item From the semi-positive $(1,1)$ form $\Theta^{(1)}$, we get $(\Theta^{(1)})^n=0$ and so the top intersection number $c_1^n(M)=0$. By Theorem \ref{main111},
$M$ can not admit a Hermitian metric with quasi-positive
Hermitian-Ricci curvature. Moreover, the quasi-positive curvature
condition in  Theorem \ref{main111} can not be replaced by
nonnegative curvature condition. \ed \eremark

\noindent Next we construct explicit Levi-Civita Ricci-flat
Hermitian metrics on all Hopf manifolds and prove Theorem
\ref{example1}.

\btheorem\label{flat} Let \beq \tilde
\omega=\omega_0-\frac{4}{n}\fR^{(1)}(\omega_0),\eeq then the first
Levi-Civita Ricci curvature of $\tilde\omega$ is zero, i.e. \beq
\fR^{(1)}(\tilde\omega)=0.\eeq \etheorem

 \bproof We consider the
perturbed metric \beq \tilde
\omega=\omega_0+4\lambda\fR^{(1)}(\omega_0),\qtq{with}
\lambda>-1,\eeq where we know form (\ref{1000001}) that
$\fR^{(1)}(\omega_0)=\sq\p\bp\log |z|^2$. That is
$$\tilde h_{i\bar j}=\frac{4}{|z|^2}\left((1+\lambda)\delta_{ij}-\frac{\lambda \bar z^i z^j}{|z|^2}\right), \qtq{and }\tilde h^{i\bar j}=\frac{|z|^2}{4}\left(\frac{\delta_{i\bar
j}}{1+\lambda}+\frac{\lambda z^i\bar
z^j}{(1+\lambda)|z|^2}\right).$$ Moreover,
$$\frac{\p \tilde h_{i\bar j}}{\p\bar z^\ell}=\frac{8\lambda z^jz^\ell\bar z^i}{|z|^6}-\frac{4(1+\lambda)\delta_{ij}z^\ell+4\lambda\delta_{i\ell} z^j}{|z|^4}\qtq{and}\frac{\p \tilde h_{i\bar j}}{\p\bar z^\ell}-\frac{\p \tilde h_{i\bar \ell}}{\p\bar z^j}=\frac{4(\delta_{i\ell}z^j-\delta_{ij}z^\ell)}{|z|^4}.$$
The Christoffel symbols of $\tilde h$ are
$$\tilde \Gamma_{\bar j i}^i=\frac{1}{2}\tilde h^{i\bar \ell}\left(\frac{\p\tilde h_{i\bar \ell}}{\p\bar z^j}-\frac{\p\tilde h_{i\bar j}}{\p\bar z^\ell}\right)
=-\frac{(n-1)z^j}{2|z|^2(1+\lambda)},$$ and $$
\p^*\tilde\omega=-2\sq\tilde \Gamma_{\bar j i}^id\bar
z^j=\sq\frac{n-1}{1+\lambda}\bp\log |z|^2, \qtq{and}
\frac{\p\p^*\tilde\omega+\bp\bp^*\tilde\omega}{2}=
\frac{n-1}{1+\lambda}\cdot\sq\p\bp\log |z|^2,$$ where the adjoint
operators $\bp^*$ and $\p^*$ are taken with respect to the new
metric $\tilde\omega$. Finally, since $\det(\tilde h_{i\bar
j})=(1+\lambda)^{n-1}4^n|z|^{-2n}$, we obtain $$ \Theta^{(1)}(\tilde
\omega)=-\sq\p\bp\log\det\tilde h_{i\bar j}=n\cdot \sq \p\bp\log
|z|^2,$$ and by Theorem \ref{main-1},
$$\fR^{(1)}(\tilde \omega)=\Theta^{(1)}(\tilde \omega)-\frac{\p\p^*\tilde\omega+\bp\bp^*\tilde\omega}{2}=\left(n-\frac{n-1}{1+\lambda}\right) \cdot \sq\p\bp\log |z|^2.$$
Now it is obvious that when $\lambda=-\frac{1}{n}$,
$\fR^{(1)}(\tilde \omega)=0$.
 \eproof
\bremark By using the same ideas and also the Ricci curvature
relations in Section \ref{riccirelations}, one can construct various
``Einstein metrics" by using different Ricci curvatures introduced
in the previous sections. \eremark

\subsection{Hermitian metrics of constant negative  scalar curvatures on Hopf
manifolds} Using the same setting as in Theorem \ref{flat}, we let
\beq \tilde \omega=\omega_0+4\lambda\fR^{(1)}(\omega_0),\qtq{with}
\lambda>-1.\eeq Hence, the Chern scalar curvature of $\tilde \omega$
is \beq \tilde s_C=\la  \Theta^{(1)}(\tilde \omega), \tilde
\omega\ra=\la n\cdot \sq \p\bp\log |z|^2,\tilde
\omega\ra=\frac{n(n-1)}{4(1+\lambda)}.\eeq On the other hand, as
computed in Theorem \ref{flat}. \beq \la
\p\p^*\tilde\omega+\bp\bp^*\tilde\omega,\tilde\omega\ra=2\left\la
\frac{n-1}{1+\lambda}\cdot\sq\p\bp\log
|z|^2,\tilde\omega\right\ra=\frac{(n-1)^2}{2(1+\lambda)^2}.\eeq
Similarly, \beq
|\p^*\tilde\omega|^2=\frac{(n-1)^2}{4(1+\lambda)^2},\eeq where the
norms are taken with respect to the new metric $\tilde \omega$.
Moreover, for the torsion $\tilde T$ of $\tilde\omega$, it is
$$\tilde T_{\ell j}^m=\tilde h^{m\bar s}\left(\frac{\p \tilde h_{j\bar s}}{\p z^\ell}-\frac{\p\tilde h_{\ell\bar s}}{\p z^j}\right)=\frac{\delta_{m\ell}\bar z^j-\delta_{mj}\bar z^\ell}{(1+\lambda)|z|^2}$$
and so \beq |\tilde T|^2=\tilde h_{m\bar n}\cdot\tilde h^{j\bar
i}\cdot\tilde h^{ \ell\bar k}\cdot\tilde T_{\ell
j}^m\cdot\bar{\tilde T^n_{ki}}=\frac{n-1}{2(1+\lambda)^2}.\eeq
Finally, by the scalar curvature relation formula (\ref{key-1}), we
see the Riemannian scalar curvature $\tilde s$ of $\tilde \omega$ is
\begin{eqnarray} \tilde s\nonumber&=&2\tilde s_C+\left(\la\p\p^*\tilde
\omega+\bp\bp^*\tilde \omega,\tilde \omega\ra-2|\p^*\tilde
\omega|^2\right)-\frac{1}{2}|\tilde
T|^2\\ \nonumber&=&2\tilde s_C-\frac{1}{2}|\tilde T|^2\\
&=&\frac{n(n-1)}{2(1+\lambda)^2}\left[\lambda-\frac{1-2n}{2n}\right].\label{scalarweight}\end{eqnarray}
It is easy to see that when $\lambda\in (-1,\infty)$, \beq \tilde
s\in (-\infty,\frac{n^2(n-1)}{4}{\Big]}.\eeq More precisely, when
$\lambda=\frac{1-n}{n}$, $\tilde s=\frac{n^2(n-1)}{4}$. Hence, for
the smooth family of Hermitian metrics $\displaystyle\tilde
\omega=\omega_0+4\lambda\fR^{(1)}(\omega_h),$ with $\lambda>-1$,
\bd\item $\lambda>\frac{1-2n}{2n}$, $\tilde\omega$ has positive
constant Riemannian scalar curvature $\tilde s$;\item when
$\lambda=\frac{1-2n}{2n}$, $\tilde\omega$ has constant zero
Riemannian scalar curvature $\tilde s$;
\item when $-1<\lambda<\frac{1-2n}{2n}$, $\tilde \omega$ has
negative constant Riemannian scalar curvature $\tilde s$. Note also
that when $\lambda\>-1$, $\tilde s\>-\infty$. \ed

\btheorem\label{example20002} For any $n\geq 2$, there exists an
$n$-dimensional compact complex manifold $X$ with $c_1(X)\geq 0$,
such that $X$ admits three different \textbf{Gauduchon metrics}
$\omega_1,\omega_2$ and $\omega_3$ with the following properties.
\bd
\item $[\omega_1]=[\omega_2]=[\omega_3]\in H^{1,1}_{A}(X)$;
\item  they  have the same
semi-positive Chern-Ricci curvature, i.e.
$$\Theta^{(1)}(\omega_1)=\Theta^{(1)}(\omega_2)=\Theta^{(1)}(\omega_3)\geq 0;$$
\item they have constant positive  Chern
scalar curvatures. \ed  Moreover,
 \bd \item $\omega_1$ has \textbf{\emph{positive}}  constant Riemannian scalar
 curvature;

 \item $\omega_2$ has \textbf{\emph{{zero}}} Riemannian scalar curvature;

\item $\omega_3$ has \textbf{\emph{negative}} constant Riemannian scalar curvature.
\ed \etheorem

\bproof Let $Y=\S^3\times \S^1$ and $\omega_0$ the canonical metric.
It is easy to see that $\p\bp\omega_0=0$, i.e. $\omega_0$ is a
Gauduchon metric. Indeed,
$$\p\bp\omega_0=\p\left(-\frac{4z^\ell}{|z|^4}d\bar z^\ell\wedge dz^i\wedge d\bar z^i\right)=\left(-\frac{4\delta_{k\ell}}{|z|^4}+\frac{8\bar z^kz^\ell}{|z|^6}\right)dz^k\wedge d\bar z^\ell\wedge dz^i\wedge d\bar z^i=0$$
since the complex dimension of $Y$ is $2$. As in the previous
paragraph, let \beq
\omega_\lambda:=\omega_0+4\lambda\fR^{(1)}(\omega_0)=\omega_0+4\sq\lambda\p\bp\log|z|^2,\qtq{with}
\lambda>-1.\eeq Hence $[\omega_\lambda]=[\omega_0]\in
H_{A}^{1,1}(Y)$ for any $\lambda>-1$, and $\omega_\lambda$ are all
Gauduchon metrics. On the other hand,
 $\omega_{\lambda}$ has first Chern-Ricci curvature
$$\Theta^{(1)}(\omega_\lambda)=2\sq \p\bp\log|z|^2\geq 0.$$
The Chern scalar curvature are \beq
s_C(\omega_\lambda)=\frac{1}{2(1+\lambda)}>0\eeq since $\lambda>-1.$
Since $n=2$, by formula (\ref{scalarweight}), the Riemannian scalar
curvature of $\omega_\lambda$ is \beq
s(\omega_\lambda)=\frac{1}{(1+\lambda)^2}\left[\lambda+\frac{3}{4}\right].\eeq
Therefore \bd\item $\lambda>-\frac{3}{4}$, $\omega_\lambda$ has
positive constant Riemannian scalar curvature;\item when
$\lambda=-\frac{3}{4}$, $\omega_\lambda$ has constant zero
Riemannian scalar curvature;
\item when $-1<\lambda<-\frac{3}{4}$, $\omega_\lambda$ has
negative constant Riemannian scalar curvature.
 \ed
Therefore, in dimension $2$, Theorem \ref{example20002} is proved.

For higher dimensional cases,  let $(Z,\omega_Z)$ be any compact
K\"ahler manifold with zero (Chern) scalar curvature. For example we
can take $(Z,\omega_Z)$ as the standard $(n-2)$-torus
$(\T^{n-2},\omega_{\T})$. Hence the background Riemannian metric has
zero Riemannian scalar curvature. Let $X$  be the product manifold
$Y\times Z$, and $\omega_X$  the product metric of $\omega_Z$ and
$\omega_\lambda$. It is obvious that $(X,\omega_X)$ is a Gauduchon
metric with semi-positive Chern-Ricci curvature and positive
constant Chern scalar curvature. Moreover, the Riemannian scalar
curvature of $\omega_X$ equals the Riemannian scalar curvature of
$\omega_\lambda$ thanks to the product structure. \eproof

\subsection{Pluriclosed metrics on the projective bundles over Hopf manifolds} Let $M=\S^{2n-1}\times \S^1$ with $n\geq
2$, and $E=T^{1,0}M$. Suppose $X=\P(E^*)$ and $L=\sO_{\P(E^*)}(1)$
is the tautological line bundle of the fiber bundle $\pi:X\>M$. By
the adjunction formula, \beq K_X=L^{-n}\ts \pi^*\left(K_M\ts \det
E\right)\eeq we see \beq K_X=L^{-n}.\eeq It is obvious, when
restricted to the fiber $X_s:=\P(E_s^*)\cong \P^{n-1}$,
$K_X|_{X_s}\cong \sO_{\P^{n-1}}(-n)$. Hence, $K_X$ is not
topologically trivial and moreover,  $c_1(X)\geq 0$. However, by a
straightforward calculation we see $c^{2n-1}_1(X)=0$ where $\dim_\C
X=2n-1$. In fact, as described in \cite{LSY12},  there is a natural
Hermitian metric $g$ on $L$ induced by the Hermitian metric
$\omega_h$ on $E=T^{1,0}M$. Since the Chern curvature tensor of
$(E=T^{1,0}M,h)$ is Griffiths-semi-positive (see formula
(\ref{Griffiths}), or \cite[Proposition~6.1]{LY12}), the curvature
$(1,1)$-form $\Theta_g$ of $(L,g)$ is semi-positive and strictly
positive on each fiber. We can see that
$\left(\Theta_{\omega_h}\right)^{2n-1}=0$ and so $c^{2n-1}_1(X)=0$.
It is easy to see that the Chern scalar curvature of $X$ is strictly
positive, i.e.
$$n\cdot \tr_g\Theta_g>0.$$
Hence, $H^0(X,mK_X)=0$ for any integer $m>0$. Now we summarize the
discussion as follows.

\bproposition Let $M=\S^{2n-1}\times \S^1$ with $n\geq 2$, and
$E=T^{1,0}M$. Suppose $X=\P(E^*)$.

\bd \item $K_X^{-1}$ admits a Hermitian metric with semi-positive
Chern curvature;

\item $K_X^{-1}$ is not topologically trivial,  i.e. $c_1(X)\neq 0$. Moreover,  $c_1(X)\geq
0$, but $c_1^{2n-1}(X)=0$;

\item $H^0(X,mK_X)=0$ for any integer $m>0$.
\ed

\eproposition

\noindent Next we consider a more general setting. Suppose $n\geq 2$
and $k\geq 1$. Let $M_n=\S^{2n-1}\times \S^1$.
$E_k=\underbrace{T^{1,0}M_n\ds\cdots \ds T^{1,0}M_n}_{\text{$k$
copies}}$ and $X_{n,k}=\P(E_k^*)\>M_n$. Since $\pi:X_{n,k}\>M_n$ is
a proper holomorphic submersion and $M_n$ can not admit any balanced
metric, by \cite[Proposition~1.9]{Mi83}, $X_{n,k}$ can not support
balanced metrics. In particular, $X_{n,k}$ is not in the Fujiki
class $\mathscr C$. On the other hand, we consider a special case
$\pi:X_{2,k}\>M_2$. It is obvious that $E_k$ has an induced
Hermitian metric with Griffiths-semi-positive curvature. Let $F$ be
the tautological line bundle of $X_{2,k}\>M_2$. The induced
Hermitian metric on $F$ has semi-positive curvature tensor
$\Theta_F$ which is also strictly positive on each fiber. Now we can
construct a family of Hermitian metrics $\omega$ with
$\p\bp\omega=0$ on $X_{2,k}$. Let $\omega_0$ be the canonical metric
(\ref{hopf}) on $M_2=\S^3\times \S^1$, and it is obvious that
$\p\bp\omega_0=0$. Then for any $\lambda>0$, \beq
\omega:=\pi^*(\omega_0)+\lambda \Theta_F\eeq  is a Hermitian metric
on $X_{2,k}$. Moreover, it satisfies $\p\bp\omega=0$.

\bproposition Suppose $n\geq 2$ and $k\geq 1$. Let
$M_n=\S^{2n-1}\times \S^1$. $E_k=\underbrace{T^{1,0}M_n\ds\cdots \ds
T^{1,0}M_n}_{\text{$k$ copies}}$ and $X_{n,k}=\P(E_k^*)$. Then

\bd \item  $X_{n,k}$ can not support any balanced metric;

\item  $X_{2,k}$ admits a Hermitian metric
$\omega$ with $\p\bp\omega=0$.

\ed

\eproposition

\section{Appendix: the Riemannian Ricci curvature and $*$-Ricci
curvature}\label{appendix}

In this appendix, we provide more details on the complexification of
Riemannian Ricci curvatures on almost Hermitian manifolds.\\

\noindent{{7.1.} \textbf{The Riemannian Ricci curvature.}}

 \blemma\label{Akey} On an almost Hermitian manifold $(M,h)$, the
 Riemannian
Ricci curvature of the background Riemannian manifold $(M,g)$
satisfies \beq Ric(X,Y)=h^{i\bar \ell}\left[R\left(\frac{\p}{\p
z^i}, X,Y, \frac{\p}{\p\bar z^\ell}\right)+R\left(\frac{\p}{\p z^i},
Y,X, \frac{\p}{\p\bar z^\ell}\right)\right]\eeq for any $X,Y\in T_\R
M$. The Riemannian scalar curvature of $(M,g)$ is \beq s=2h^{i\bar
j}h^{k\bar \ell}\left(2R_{i\bar \ell k\bar j}-R_{i\bar j k\bar
\ell}\right).\eeq \elemma

\bproof For any $X,Y\in T_{\R}M$, by using real coordinates $\{x^i,
x^I\}$ and the relations (\ref{1001}), (\ref{1002}), (\ref{1003})
and (\ref{1004}), one can see: \be &&Ric(X,Y)\\&=&
g^{i\ell}R\left(\frac{\p}{\p x^i}, X,Y, \frac{\p}{\p x^\ell}\right)+g^{iL}R\left(\frac{\p}{\p x^i}, X,Y, \frac{\p}{\p x^L}\right)\\&&+g^{I\ell}R\left(\frac{\p}{\p x^I}, X,Y, \frac{\p}{\p x^\ell}\right)+g^{IL}R\left(\frac{\p}{\p x^I}, X,Y, \frac{\p}{\p x^L}\right)\\
&=&g^{i\ell}R\left(\frac{\p}{\p z^i}+\frac{\p}{\p\bar z^i}, X,Y,
\frac{\p}{\p z^\ell}+\frac{\p}{\p\bar z^\ell}\right)+\sq
g^{iL}R\left(\frac{\p}{\p z^i}+\frac{\p}{\p\bar z^i}, X,Y,
\frac{\p}{\p z^\ell}-\frac{\p}{\p\bar z^\ell}\right)\\
&&+\sq g^{I\ell}R\left(\frac{\p}{\p z^i}-\frac{\p}{\p\bar z^i}, X,Y,
\frac{\p}{\p z^\ell}+\frac{\p}{\p\bar
z^\ell}\right)-g^{IL}R\left(\frac{\p}{\p z^i}-\frac{\p}{\p\bar z^i},
X,Y,
\frac{\p}{\p z^\ell}-\frac{\p}{\p\bar z^\ell}\right)\\
&=&2g^{i\ell}\left[R\left(\frac{\p}{\p z^i}, X,Y, \frac{\p}{\p\bar
z^\ell}\right)+R\left(\frac{\p}{\p \bar z^i}, X,Y, \frac{\p}{\p
z^\ell}\right)\right]\\
&&+2\sq g^{iL}\left[-R\left(\frac{\p}{\p z^i}, X,Y, \frac{\p}{\p\bar
z^\ell}\right)+R\left(\frac{\p}{\p \bar z^i}, X,Y, \frac{\p}{\p
z^\ell}\right)\right]\\
&=&h^{i\bar \ell}R\left(\frac{\p}{\p z^i}, X,Y, \frac{\p}{\p\bar
z^\ell}\right)+h^{\ell\bar i}R\left(\frac{\p}{\p \bar z^i}, X,Y,
\frac{\p}{\p z^\ell}\right)\\
&=&h^{i\bar \ell}\left[R\left(\frac{\p}{\p z^i}, X,Y,
\frac{\p}{\p\bar z^\ell}\right)+R\left(\frac{\p}{\p z^i}, Y,X,
\frac{\p}{\p\bar z^\ell}\right)\right].\ee

\noindent By using the symmetry $R(X,Y,Z,W)=R(Y,X,W,Z)$ for
$X,Y,Z,W\in T_{\C} M$, we get the following formula for the
{Riemannian scalar curvature}, \be s &=&h^{i\bar j}h^{k\bar
\ell}\left( R_{ik\bar \ell\bar j}+ R_{i\bar \ell k\bar j}+ R_{\bar j
k\bar \ell i}+ R_{\bar j \bar \ell k i}\right)\\&=&2h^{i\bar
j}h^{k\bar \ell}\left( R_{ik\bar \ell\bar j}+ R_{i\bar \ell k\bar
j}\right)\\&=&2h^{i\bar j}h^{k\bar \ell}\left(2R_{i\bar \ell k\bar
j}-R_{i\bar j k\bar \ell}\right).\ee The proof of Lemma \ref{Akey}
is complete. \eproof

{\noindent {7.2.} \textbf{$*$-Ricci curvature and $*$-scalar
curvature.}} \bdefinition Let $\{e_i\}_{i=1}^{2n}$ be an orthonormal
basis of $(T_\R M,g)$, the (real) $*$-Ricci curvature of $(M,g)$ is
defined to be (e.g. \cite{TV81}) \beq
Ric^*(X,Y):=\sum_{i=1}^{2n}R(e_i,X,JY,Je_i), \eeq for any $X,Y\in
T_\R M$. The $*$-scalar curvature (with respect to the Riemannian
metric) is defined to be \beq s^*=\sum_{j=1}^{2n}Ric^*(e_i,e_i).\eeq
\edefinition

\noindent It is easy to see that, for any $X,Y\in T_\R M$.
 \beq
Ric^*(X,Y)=Ric^*(JY,JX).\label{skew*Ric}\eeq \blemma We have the
following formula for $*$-Ricci curvature,
 \begin{eqnarray}
Ric^*(X,Y)\label{*Ric3}&=&\sq h^{\ell\bar i} R\left(\frac{\p}{\p\bar
z^i},X,JY,\frac{\p}{\p z^\ell}\right)-\sq h^{i\bar \ell}
R\left(\frac{\p}{\p z^i},X,JY,\frac{\p}{\p \bar
z^\ell}\right)\\
&=& \sq h^{k\bar \ell} R\left(\frac{\p}{\p z^k}, \frac{\p}{\p\bar
z^\ell},X,JY\right).\label{*Ric}\end{eqnarray}

\noindent The $*$-scalar curvature is \beq s^*=2h^{i\bar j}h^{k\bar
\ell}R_{i\bar j k\bar \ell}=2s_H.\eeq\elemma

\bproof As similar as the computations in Lemma \ref{Akey}, we can
write down the $*$-Ricci curvature in real coordinates
$\{\frac{\p}{\p x^i}, \frac{\p}{\p x^I}\}$ and show: \be
Ric^*(X,Y)&=&
\sq h^{\ell\bar i}R\left(\frac{\p}{\p\bar z^i},X,JY,\frac{\p}{\p
z^\ell}\right)-\sq h^{i\bar \ell} R\left(\frac{\p}{\p
z^i},X,JY,\frac{\p}{\p \bar z^\ell}\right).\ee

\noindent  On the other hand, \be &&Ric^*(X,Y)+Ric^*(JY,JX)\\&=&\sq
h^{\ell\bar i} R\left(\frac{\p}{\p\bar z^i},X,JY,\frac{\p}{\p
z^\ell}\right)-\sq h^{i\bar \ell} R\left(\frac{\p}{\p
z^i},X,JY,\frac{\p}{\p \bar z^\ell}\right)\\
&&-\sq h^{\ell\bar i} R\left(\frac{\p}{\p\bar z^i},JY,X,\frac{\p}{\p
z^\ell}\right)+\sq h^{i\bar \ell}
R\left(\frac{\p}{\p z^i},JY,X,\frac{\p}{\p \bar z^\ell}\right)\\
&=&2\sq h^{k\bar \ell} R\left(\frac{\p}{\p z^k}, \frac{\p}{\p\bar
z^\ell},X,JY\right), \ee where the last step follows by Bianchi
identity. Therefore (\ref{*Ric}) follows by the relation
(\ref{skew*Ric}). For the scalar curvature $s^*$, by definition, it
is \be s^*&=&g^{jk}Ric^*\left(\frac{\p}{\p x^j}, \frac{\p}{\p
x^k}\right)+g^{jK}Ric^*\left(\frac{\p}{\p x^j}, \frac{\p}{\p
x^K}\right)\\&&+g^{Jk}Ric^*\left(\frac{\p}{\p x^J}, \frac{\p}{\p
x^k}\right)+g^{JK}Ric^*\left(\frac{\p}{\p x^J}, \frac{\p}{\p
x^K}\right).\ee By using the symmetry $R(X,Y,Z,W)=R(Y,X,W,Z)$ and
(\ref{*Ric3}), we see \be s^*&=&(\sq h^{\ell \bar i})\left[\sq
h^{k\bar j}R\left(\frac{\p}{\p \bar z^i},\frac{\p}{\p\bar
z^j},\frac{\p}{\p z^k},\frac{\p}{\p z^\ell}\right)-\sq h^{j\bar
k}R\left(\frac{\p}{\p \bar z^i},\frac{\p}{\p z^j},\frac{\p}{\p \bar
z^k},\frac{\p}{\p
z^\ell}\right)\right]\\
&&-(\sq h^{i\bar \ell})\left[\sq h^{k\bar j}R\left(\frac{\p}{\p
 z^i},\frac{\p}{\p\bar z^j},\frac{\p}{\p z^k},\frac{\p}{\p
\bar z^\ell}\right)-\sq h^{j\bar k}R\left(\frac{\p}{\p
z^i},\frac{\p}{\p z^j},\frac{\p}{\p \bar z^k},\frac{\p}{\p \bar
z^\ell}\right)\right]\\
&=&-h^{i\bar j}h^{k\bar \ell}R_{\bar j \bar \ell k i}+h^{i\bar
j}h^{k\bar \ell}R_{\bar j k\bar\ell i}+h^{i\bar j}h^{k\bar \ell}
R_{i\bar\ell k\bar j}-h^{i\bar j}h^{k\bar\ell}R_{ik\bar\ell\bar
j}\\
&=&2h^{i\bar j}h^{k\bar \ell}R_{i\bar j k\bar \ell}=2s_H,\ee where
the last step follows by Bianchi identity.
 \eproof

\bremark\label{remark}  It is easy to see that $*$-Ricci curvature
is neither symmetric nor skew-symmetric. For example, by using the
Hermitian Ricci tensor and (\ref{*Ric}), we see the following
submatrix of the real matrix representation of $*$-Ricci curvature:
\begin{eqnarray} Ric^*\left(\frac{\p}{\p x^i}, \frac{\p}{\p
x^j}\right)\nonumber&=& \sq h^{k\bar \ell} R\left(\frac{\p}{\p z^k},
\frac{\p}{\p\bar z^\ell},\frac{\p}{\p x^i}, \frac{\p}{\p
x^J}\right)\\ \nonumber &=&- h^{k\bar \ell} R\left(\frac{\p}{\p
z^k}, \frac{\p}{\p\bar z^\ell},\frac{\p}{\p z^i}+\frac{\p}{\p\bar
z^i}, \frac{\p}{\p z^j}-\frac{\p}{\p\bar z^j}\right)\\ \nonumber
&=&\left(R_{\bar i\bar j}-R_{ij}\right)+(R_{i\bar j}-R_{\bar i
 j})\\ &=&\left(R_{\bar i\bar j}-R_{ij}\right)+(R_{i\bar j}+R_{j\bar i
 }).\label{*Ric2}\end{eqnarray} The first part in (\ref{*Ric2}) is skew
symmetric whereas the second part is symmetric. Hence, as a real
$(0,2)$ tensor,  it is impossible to define the positivity or
negativity for the $*$-Ricci curvature. That is $Ric^*(X,X)>0$ for
all nonzero vector $X\in T_\R M$ can not happen on any almost
Hermitian manifold.  To make the $*$-Ricci tensor a symmetric
tensor, an extra condition  as \beq R(X,Y,Z,W)=R(X,Y,JZ,JW)\eeq is
sufficient (see, e.g. \cite{Gr76,TV81}).

 \eremark

\end{document}